\documentclass[final,onefignum,onetabnum]{siamonline171218}

\usepackage{amsmath,bm}

\usepackage{amsfonts}
\usepackage{epstopdf}
\usepackage{graphicx}
\usepackage{url}
\usepackage{dsfont}
\usepackage{xcolor}
\usepackage{subcaption}
\usepackage{setspace}
\usepackage{booktabs}

\DeclareMathAlphabet\mathpzc{OT1}{pzc}{m}{it}
\let\mathcal=\mathpzc

\def\P{{\mathbb P}}

\def\pr{{\mathrm{pr}}}
\def\erf{{\mathrm{erf}}}

\let\<=\langle
\let\>=\rangle
\let\^=\hat
\def\~#1{{\mbox{#1}}}
\def\N{{\mathbb N}}
\def\R{{\mathbb R}}
\let\-=\mathbf

\def\circ{\ifmmode\mathchar"220E\else$\mathchar"220E$\fi}

\def\@#1{{\bf #1}}

\usepackage{amsfonts}
\usepackage{graphicx}
\usepackage{epstopdf}
\usepackage{algorithm}

\usepackage{algpseudocode,float}

\ifpdf
  \DeclareGraphicsExtensions{.eps,.pdf,.png,.jpg}
\else
  \DeclareGraphicsExtensions{.eps}
\fi

\usepackage{enumitem}
\setlist[enumerate]{leftmargin=.5in}
\setlist[itemize]{leftmargin=.5in}

\newsiamremark{remark}{Remark}
\newsiamremark{hypothesis}{Hypothesis}
\crefname{hypothesis}{Hypothesis}{Hypotheses}
\newsiamthm{claim}{Claim}

\headers{Bayesian inference for Poisson data}{Q. Zhou, T. Yu, X. Zhang and J. Li}

\title{Bayesian inference and uncertainty quantification for medical image reconstruction
with Poisson data\thanks{Submitted to the editors DATE.
\funding{This work was funded by the National Natural Science Foundation of China, under grant numbers 11771288 and 11771289.
JL also wants to acknowledge the support of the Student Innovation Center at Shanghai Jiao Tong University.}}}

\author{Qingping Zhou\thanks{School of Mathematical Sciences, Institute of Natural Sciences,
Shanghai Jiao Tong University, 800 Dongchuan Rd, Shanghai 200240, China, (zhou2015@sjtu.edu.cn).}
\and Tengchao Yu\thanks{School of Mathematical Sciences, Institute of Natural Sciences,
Shanghai Jiao Tong University, 800 Dongchuan Rd, Shanghai 200240, China, (tengchaoyu@sjtu.edu.cn).}
\and Xiaoqun Zhang\thanks{Institute of Natural Sciences, School of Mathematical Sciences, and
the MOE Key Laboratory of Scientific and Engineering Computing,
Shanghai Jiao Tong University, 800 Dongchuan Rd, Shanghai 200240, China, (xqzhang@sjtu.edu.cn)..}
\and Jinglai Li\thanks{Corresponding Author, Department of Mathematical Sciences, University of Liverpool,
Liverpool, L69 7ZL, UK (jinglai.li@liverpool.ac.uk).}}

\usepackage{amsopn}

\ifpdf
\hypersetup{
  pdftitle={An Example Article},
  pdfauthor={Q. Zhou, X. Zhang and J. Li}
}
\fi

\begin{document}

\maketitle

\begin{abstract}
  We provide a complete framework for performing infinite-dimensional Bayesian inference and uncertainty quantification for image reconstruction with Poisson data. In particular, we address the following issues to make the Bayesian framework applicable in practice. 
We first introduce a positivity-preserving reparametrization, and we prove that  under the reparametrization
and a hybrid prior, the posterior distribution is well-posed in the infinite dimensional setting.   
Second we provide a dimension-independent MCMC algorithm, based on  the preconditioned Crank-Nicolson Langevin  method, in which we use a primal-dual scheme to compute the offset direction. 
Third we give a method combining the model discrepancy method
and maximum likelihood estimation to determine the regularization parameter in the hybrid prior.
Finally we propose to use the obtained posterior distribution to detect artifacts in a recovered image. 
We provide an example to demonstrate the effectiveness of the proposed method.  
\end{abstract}



\section{Introduction}\label{s:intro}
Image reconstruction involves constructing interpretable images of objects of interest from the data recorded by an imaging
device~\cite{fessler2017medical}. Image reconstruction is usually cast as an inverse problem as one wants to determine the input to a
system from the output of it. 
In most practical image reconstruction problems, the measurement and recording process is inevitably corrupted by noise,
which renders the obtained data random.  
The statistical properties of the data have significant impact to the reconstruction results. 
In this work we shall focus on a special type of medical image reconstruction problems where the recorded data follows a Poisson distribution. 
The Poisson data usually arises in imaging problems where the unknown quantity of interest is an object which interacts with some known incident beam of
photons or electrons~\cite{hohage2016inverse}. 
A very important example of such problems is the Positron emission tomography(PET)~\cite{ollinger1997positron,bailey2005positron},  a nuclear medicine imaging technique that is widely used in early detection and treatment follow up of many diseases, including cancer.
In PET, the detection of signal is essential a photon counting process and as a result the data is well modeled by a Poisson distribution~\cite{bailey2005positron,hohage2016inverse}. 
The problem has attracted considerable research interests, and a number of methods have been developed to recover the image, e.g., \cite{shepp1982maximum,vardi1985statistical,fessler1994penalized}, just to name a few.

On the other hand, the stochastic nature of the data also introduces uncertainty into the image reconstruction process, and as a result the image obtained is unavoidably subject to uncertainty.
In practice, many important decisions such as diagnostics have to be made based on the images obtained. 
It is thus highly desirable to have methods that can not only compute the image but also quantify the uncertainty in the image obtained. 
To this end, the Bayesian inference method has become a popular tool for image reconstruction \cite{kaipio2005statistical}, largely thanks to its ability to quantify uncertainty 
in the obtained image. 
The Bayesian formulation has long been used to solve image reconstruction problems with Poisson data, e.g.,~\cite{hebert1989generalized,mumcuoglu1996bayesian,green1990bayesian}.
We note, however, that most of the works in the early years focus on computing a point estimate, which is usually the maximum a posteriori (MAP) estimate in the Bayesian setting, because of the limited computational power available then. 
More recently, mounting interest has been directed to the computation of the complete posterior distribution, rather than a point estimate,
of the image, for that it can provide the important uncertainty information of the reconstruction results.  
For example, a Markov chain Monte Carlo (MCMC) algorithm is developed to sample the posterior distribution of the image in~\cite{bardsley2016metropolis}, and a variational Gaussian approximation of the posterior is proposed in~\cite{arridge2018variational}. 

A serious challenge in the numerical implementation of the Bayesian image reconstruction is that  in certain circumstances  the inference results diverge with respect to resolution/discretization  refinement, which is known as to be
discretization variant or dimension dependent. 
To address the issue, Stuart~\cite{stuart2010inverse} proposes an infinite dimensional  framework, formulating the Bayesian inference problem in function spaces.
{Under the infinite dimensional framework, the inference results will converge with respect to discretization dimensionality, which is an important property for the numerical implementation. 
For example, it allows one to use multigrid strategy, e.g.~\cite{ye2001nonlinear,pan1991numerical}, to accelerate the sampling of the posterior. }
Building on several existing works, we aim to provide in this work a complete framework for performing infinite dimensional Bayesian inference and uncertainty quantification 
for medical image reconstruction with Poisson data, while providing treatments of  several issues surrounding the problem. Specifically we summarize the key ingredients of our Bayesian framework as the following. 
First, in the usual setup, the function of interest can be both positive and negative valued. 
However, in the Poisson problem, when the function is negative valued, it may cause the Poisson likelihood function to be undefined (see Section~\ref{sec:prior} for more details),
which renders the posterior distribution ill-posed in the infinite dimensional setting. 
To tackle the issue,  we introduce a reparametrization of the unknown image which ensures that the function of interest is always positive valued. 
Moreover, medical images are often subject to sharp jumps and here
 we use the TV-Gaussian (TG) hybrid prior distribution proposed in \cite{yao2016tv} to model the jumps in the function/image. 
Using the positivity-preserving reparametrization and the TG prior, we are able to show that the resulting posterior distribution is \emph{well-posed} in the infinite dimensional setting, which, to the best of our knowledge, has not yet been done for the Poisson data model.  
Second, we consider the numerical implementation of the Bayesian inference. 
A main difficulty here is that 
many standard MCMC algorithms such as the well-known Metropolis-Hastings (MH)~\cite{robert2013monte,brooks2011handbook}, degenerate with respect to resolution refinement. 
In~\cite{cotter2013mcmc} the authors introduce a MCMC algorithm termed as the preconditioned Crank-Nicolson (pCN) method, 
the performance of which is independent of discretization dimensionality. 
The authors also provide a Langevin variant of the pCN algorithm in \cite{cotter2013mcmc} which accelerates the sampling procedure
by incorporating the local gradient information of the likelihood function. 
In our problem, the pCN-Langevin (pCNL) algorithm can not be used directly because the prior used here has the total variation (TV) term which can be  non-differentiable. 
To overcome this difficulty we modify the pCNL method by replacing the gradient direction with one computed by the primal-dual algorithm.  
We note that a similar problem is considered in \cite{pereyra2016proximal,durmus2018efficient} where a proximal method is used to approximate the gradient direction. Other than that the directions are computed with different approaches, 
another main difference between the aforementioned works and the present one is that we use the pCN framework here so the algorithm is dimension independent, while the works \cite{pereyra2016proximal,durmus2018efficient} concern finite dimensional problems where discretization refinement is not an issue. 
Third, an important issue in the TG hybrid prior is to determine the value of the regularization parameter of the TV term. 
In the Bayesian framework, such parameters are often determined with the hierarchical Bayes or the empirical Bayes method~\cite{gelman2013bayesian}. 
As discussed in Section~\ref{sec:lambda}, these methods, however, are computationally intractable in our problem as we do not know the normalization constant of the TG prior. 
{Thus, in this work we provide a method to determine the value of the TV regularization parameter by combining the realized discrepancy
model fit assessment approach developed in \cite{Gelman1996Posterior}
and the stochastic proximal gradient method developed in \cite{fernandez2018maximum}.}
Finally, we provide an application of the uncertainty information obtained in the Bayesian framework,
where we use the posterior distribution to detect possible artifacts in any reconstructed image.

The rest of the paper is organized as the following.
In Section~\ref{sec:formulation}, we present the infinite dimensional Bayesian formulation of the image reconstruction problem
 with Poisson data, and we prove that under the reparametrization the resulting posterior is well-posed in the function space.   
In Section~\ref{sec:mcmc}, we describe the primal-dual pCN algorithm to sample the posterior distribution of the present problem.
Section~\ref{sec:lambda} provides a method to determine the value of the regularization parameter  in the TG prior. 
Section~\ref{sec:artifact} discusses how to use the posterior distribution to detect artifacts in a reconstructed image. 
Finally numerical experiments of the proposed Bayesian framework are performed in Section~\ref{sec:results}.

  \section{Infinite dimensional Bayesian image reconstruction with Poisson data}\label{sec:formulation}
In this section, we formulate the image reconstruction with Poisson data in an infinite dimensional Bayesian framework.

\subsection{The Bayesian inference formulation for functions}
We start by presenting a generic Bayesian inference problem for functions.
Let $X$ be a separable {Hilbert} space of functions with inner product $\<\cdot,\cdot\>_X$.
Our goal is to infer $u\in X$ from data $y\in Y\subset R^{d}$  and $\-y$ is related to $u$ via the likelihood function $\pi(\-y|u)$,
 i.e., the distribution of $\-y$ conditional on the value of $u$.
In the Bayesian setting,
 we first assume a prior distribution $\mu_\mathrm{pr}$ of the unknown $u$, which represents one's prior knowledge on the unknown.
 In principle $\mu_\mathrm{pr}$ can be any probabilistic measure defined on the space $X$.
The posterior measure $\mu^y$ of $u$ conditional on data $y$
is provided by the Radon-Nikodym(R-N) derivative:
\begin{equation} \frac{d\mu^y}{d\mu_\mathrm{pr}}(u) = \pi(\-y|u), \label{e:bayes}
\end{equation}
which can be interpreted as the Bayes' rule in the infinite dimensional setting.
The posterior distribution $\mu^y$ thus can be computed from Eq.~\eqref{e:bayes} with, for example, a MCMC simulation.

\subsection{Poisson data model and the positivity-preserving reparametrization}\label{s:setup}
To perform the Bayesian inference, we first need to specify the likelihood function, which can be derived from the underlying mathematical model relating the data and to the unknown image.
%
We assume that the image is first projected to the noise-free observable via a mapping $A: X \rightarrow Y$,
\begin{subequations}\label{e:radon}
\begin{equation}
\bm{\theta} = \~A u.
\end{equation}
While noting that the proposed framework is rather general, here for simplicity we restrict ourselves in the cases where $A$ is a bounded linear transform. For example in the PET imaging problems, the mapping $A$ is {approximately} the  Radon transform,
where each $\theta_i$ is computed by integrating $u(\-x)$ alone a line $L_i$:
\begin{equation}
\theta_i= (\~Au)_i = K\int_{L_i} u(\-x) |d\-x|,
\end{equation}
\end{subequations}
for $i = 1...d$, where $K$ is a positive constant describing the noise level.
Poisson noise is then applied to the \textit{projected observable} $\bm{\theta}$, yielding the likelihood function $\pi(\-y|u)= \pi_\mathrm{P}(\-y|\bm{\theta}=\~A u)$,
where $\pi_\mathrm{P}(\-y|\bm\theta)$ is the $d$-dimensional Poisson distribution:
\begin{equation}
 \pi_\mathrm{P}(\-y|\bm\theta) = \prod_{i=1}^{d} \frac{ (\theta_i)^{y_i} \exp(-\theta_i)}{y_i!}. \label{e:ytheta}
\end{equation}
In the PET problem, there is an additional restriction: the unknown function $u$ must be positive.
The reason is two-fold: first from the physical point of view, the unknown $u$ represents the density of the medium, which is positive;
from a technical point of view, if $u$ is not constrained to be positive, it may yield some negative components of the predicted data $\bm{\theta}$, which renders the Poisson likelihood un-defined.
To this end, we need to introduce a transformation to preserve positivity of the unknown $u$.
To impose the positivity constraint, we reparameterize the  unknown $u$ as: 
\[u(\-x) = f(z(\-x)) =\frac{a}2(\erf(z(\-x))+b),\]
where $a$ and $b$ are two constants satisfying $a>0$ and $b>1$, and $\erf(\cdot)$ is the error function defined as:
\[u(\-x) = f(z(\-x)) =\frac{a}2(\erf(z(\-x)/c)+b),\]
where $a$ $b$ and $c$ are constants satisfying $a>0$, $b>1$, and $c>0$, and $\erf(\cdot)$ is the error function defined as:
 \[
\erf(z) = \frac{2}{\sqrt{\pi}}\int_0^z e^{-t^2} dt.
\]
With the new parametrization, it is easy to see that for any $\-x\in \Omega$, we have
\begin{equation}
a(b-1)\leq u(\-x)\leq a(b+1). \label{e:ubound}
\end{equation}
{Moreover, as the behavior of the error function is well understood~(for example
its derivative is simply the Gaussian distribution),
which allows us determine the parameters conveniently.
That said, it is worth noting here that the methods presented here
does not rely on this specific reparametrization formulation.}
Now we can infer the new unknown $z$ and once $z$ is known $u$ can be computed accordingly. 
In this setup, the likelihood function for $z$ becomes
\[ 
\pi(\-y|z) = \pi_P(\-y|\bm\theta= \~A f(z))
\]
where $\pi_P(\-y|\bm\theta)$ is the $d$-dimensional Poisson distribution given by Eq.~\eqref{e:ytheta}. 
Following \cite{stuart2010inverse}, we can write the likelihood function $\pi(\-y|z)$ in the form of
\begin{subequations}\label{e:lh}
\begin{equation}
\pi(\-y|z) \propto \exp(-\Phi(z,\-y)),
\end{equation}
where
\begin{equation}
\Phi(z;\-y) = \sum_{i=1}^d (\~A f(z))_i - y_i \ln(\~A f(z))_i. \label{e:phiz}
\end{equation}
\end{subequations}
For simplicity we can rewrite Eq.~\eqref{e:phiz} as,
\begin{equation} \Phi(z;\-y) = \langle \~Af(z), \bm{1}\rangle - \langle \-y,\ln(\~Af(z)) \rangle
= \langle \bm{\theta}, \bm{1}\rangle - \langle \-y,\ln \bm{\theta} \rangle, \label{e:phiz2} \end{equation}
where $\bm{1}$ is a $d$-dimensional vector whose components are all one, $\bm{\theta}=\~Af(z)$ is the predicted observable,
and $\<\cdot,\cdot\>$ denotes the Euclidean inner product.
In what follows we often omit the argument $y$ and simply use $\Phi(z)$, when not causing ambiguity.
This notation will be used often later.


\subsection{Bayesian framework with the hybrid prior for PET imaging}
\label{sec:prior}
We now describe how the infinite dimensional Bayesian inference framework is applied to the PET problem.
First we assume that the unknown function $z$ is a function defined on $\Omega$, a bounded open subset of $R^2$.
In particular, we set the state space $X$ to be the Sobolev space $H^1(\Omega)$:
\[X=H^1(\Omega) = \{ z(\-x) \in L_2(\Omega)\, |\, {\partial_{x_1} z},\,{\partial_{x_2} z} \in L_2(\Omega) \,\mathrm{for\,all}\, \-x=(x_1,\,x_2)\in \Omega\, \}.\]
The associated norm $\|\cdot\|_X=\|\cdot\|_{H^1}$ is
\[ \|z\|_{H^1} ^2= \|z\|^2_{L_2(\Omega)}+\|\partial_{x_1}z\|^2_{L_2(\Omega)}+\|\partial_{x_2}z\|^2_{L_2(\Omega)}.\]
Choosing a good prior distribution is one of the most important issues in
 Bayesian inference.
Conventionally one often assumes that the prior on $z$, is a  Gaussian measure defined on $X$ with mean $\xi$ covariance operator $C_0$,
i.e. {$\mu_\mathrm{pr}= N(\xi,C_0)$}.
Note that $C_0$ is symmetric positive and of trace class.
The Gaussian prior has many theoretical and practical advantages, but a major limitation of the Gaussian prior is that it can not model functions with sharp jumps well.

To address the issue, here we use the TV-Gaussian prior proposed in \cite{yao2016tv}:
\begin{equation}
\frac{d\mu_\pr}{d\mu_{0}}(z) \propto\exp(-R(z)),\quad R(z) = \lambda\| z\|_\textsc{tv}.\label{e:priorz}
\end{equation}
where {$\mu_0= N(\xi,C_0)$ is the Gaussian reference prior defined on $X$ with mean $\xi$ and covariance $C_0$} and $\|\cdot\|_\mathrm{TV}$ is the TV seminorm,
\begin{equation}
\| z\|_\mathrm{TV} = \int_\Omega \|\nabla u\|_2 dx, \label{e:tv}    
\end{equation}
and $\lambda$ is a positive constant.
It follows immediately that the R-N derivative of $\mu^y$ with respect to $\mu_0$ is
\begin{equation}
\frac{d\mu^y}{d\mu_0}(z) \propto \exp(-\Phi(z)-R(z)) , \label{e:bayes2}
\end{equation}
which returns to the conventional formulation of inference with a Gaussian prior.
Thus all the methods developed for inference problems with Gaussian priors can be directly applied to our formulation.
{ We note that
it is natural to directly apply the TV seminorm  to the original image $u$; 
if we do so, however, Proposition~\ref{prop:phi} may no longer hold. 
For this technical reason we here choose to impose the TV seminorm on the new variable $z$. 
Nonetheless, it can be showed that 
\[
\|z\|_\mathrm{TV} \geq
\frac{a}c \|u\|_\mathrm{TV}.
\]
}

Next we shall show that the formulated Bayesian inference problem is well defined in the infinite dimensional setting.
We first show that $\Phi(z)$ given by Eq.~\eqref{e:phiz2} satisfies certain important conditions,  as is stated by Proposition \ref{prop:phi}.
\begin{proposition} \label{prop:phi}
The functional $\Phi$ given in Eq.~\eqref{e:phiz2} has the following properties:
\begin{enumerate}
\item For every $r>0$, there are  constants $ M(r) \in \R$ and $N(r) > 0$ such that, for all $z\in X$ , and $y \in Y$ with $\|y\|_2<r$,
\[M \leq \Phi(z)\leq N;\]
\item For every $r>0$ there is  a constant $M(r)>0$ such that, for all $z,\,v\in X$ with \\
$\max\{\|z\|_X,\|v\|_X\}<r$,
\[|\Phi(z)-\Phi(v)|\leq M \|z-v\|_X;\]
 \item There exists a constant $M>0$ such that for any $y,\,y'\in Y$,  we have
    \[     |\Phi(z;y)-\Phi(z;y')|\leq M \|y-y'\|_2.   \]
\end{enumerate}
\end{proposition}

A detailed proof of the proposition is provided in Appendix A.
Following {Proposition} \ref{prop:phi}, we can conclude that the hybrid prior~\eqref{e:priorz}, and the log-likelihood function Eq.~\eqref{e:phiz2} yield a well-behaved posterior measure given by Eq.~\eqref{e:bayes2} in the infinite-dimensional setting, as is summarized in the following theorem:
\begin{theorem}\label{th:wellposed}
For $\Phi(z)$ given by Eq.~\eqref{e:phiz2} and prior measure $\mu_\mathrm{pr}$ given by Eq.~\eqref{e:priorz}, we have the following results:
\begin{enumerate}
\item $\mu^y$ given by Eq.~\eqref{e:bayes2} is a well-defined probability measure on $X$. \label{th:1}
\item $\mu^y$ given by Eq.~\eqref{e:bayes2} is Lipschitz in the data $y$, with respect to the Hellinger distance: if $\mu^y$ and $\mu^{y'}$ are two measures corresponding to data $y$ and $y'$ the there exists $C=C(r)$ such that, for all $y,y'$ with $\max\lbrace\Vert y\Vert_2,\Vert y'\Vert_2\rbrace<r$,
\[
d_\mathrm{Hell}(\mu^y,\mu^{y'})\leq C\Vert y-y'\Vert_2.
\]
\item
Let
\begin{equation}
\frac{d\mu_{N_1,N_2}^y}{d\mu_0} = \exp (-\Phi_{N_1}(z)-R_{N_2}(z)),
\end{equation}
where $\Phi_{N_1}(z)$ is a $N_1\in \N$ dimensional approximation of $\Phi(z)$ and $R_{N_2}(z)$ is a $N_2\in\N$ dimensional approximation of $R(z)$.
Assume that $\Phi_{N_1}$ satisfies the three properties of Proposition~\ref{prop:phi} with constants uniform in $N_1$, and  $R_{N_2}$ satisfy Assumptions A.2 (\romannumeral1) and (\romannumeral2) in \cite{yao2016tv}
 with constants uniform in $N_2$. Assume also that for {any} $\epsilon>0$, there exist two positive sequences $\{a_{N_1}(\epsilon)\}$ and $\{b_{N_2}(\epsilon)\}$ converging to zero, such that
$\mu_0(X_\epsilon)\geq1-\epsilon$ for {any} $N_1,N_2\in \N$,
where
\[X_\epsilon 
= \{z\in X\, |\,\vert\Phi(z)-\Phi_{N_1}(z)\vert\leq a_{N_1}(\epsilon),\,
\vert R(z)-R_{N_2}(z)\vert\leq b_{N_2}(\epsilon) \}.
\]
Then we have
$$
d_\mathrm{Hell}(\mu^y,\mu^y_{N_1,N_2})\to0 ~~~~\mathrm{as}~~~~ N_1,\,N_2\to+\infty.
$$
\end{enumerate}
\end{theorem}
Theorem~\ref{th:wellposed} is a direct consequence of Proposition~\ref{prop:phi}, and the proof of the theorem can be found in \cite{yao2016tv} and is omitted here.

Finally, we note that, in the numerical implementations,
we use the truncated Karhunen-Lo$\grave{\textrm{e}}$ve (KL) expansion~\cite{Li2015note}  to represent the unknown $z$.
Namely, we write $z$ as
\begin{equation}
z(\-x) = \sum_{i=1}^N z_i \sqrt{\eta_i} e_i(\-x),
\end{equation}
where $\{\eta_i,\,e_i(\-x)\}$ are the eigenvalue-eigenfunction pair of the covariance operator $C_0$,
and $(z_1, \,...,\, z_N)$ are independent  with each following a standard normal distribution.
In the KL representation, the number of KL modes (eigenfunctions) $N$ corresponds to the discretization dimensionality.

\section{The primal-dual preconditioned Crank-Nicolson MCMC algorithm}\label{sec:mcmc}
In most practical image reconstruction problems, the posterior distribution can not be analytically calculated. Instead, one usually represent the posterior
by samples drawn from it using MCMC algorithms.
It is demonstrated in~\cite{cotter2013mcmc} that standard MCMC algorithms may become problematic in the infinite dimensional setting: its acceptance probability will generate to zero as the discretization dimensionality increases.
Here we adopt the pCN MCMC algorithm particularly developed for the infinite dimensional problems~\cite{cotter2013mcmc}.
An important feature of the pCN MCMC algorithm is that its sampling efficiency is independent of discretization dimensionality up to the numerical errors in the evaluation of the functionals $R(\cdot)$ and $\Phi(\cdot)$, which makes it particular useful for sampling the posterior distribution defined in function spaces.
We start with a brief introduction of the pCN algorithm following the presentation of \cite{cotter2013mcmc}.
We denote $\Phi(z)+R(z)$ of Eq. \eqref{e:bayes2} as $\Psi(z)$.
Simply speaking the algorithms are derived by applying the Crank-Nicolson (CN) discretization to a stochastic partial differential equation whose
invariant distribution is the posterior.
We here omit the derivation details while referring interested readers to \cite{cotter2013mcmc}, and jump directly to the pCN proposal:
\begin{equation}
v=(1-\beta^2)^{\frac12}z+ \beta w, \label{e:pcn}
\end{equation}
where $z$ and $v$ are the present and the proposed positions respectively,  $w\sim N(\xi,C_0)$ and $\beta\in[0,1]$ is the parameter controlling the stepsize of the algorithm.
The proposed sample $v$ is then accepted or rejected according to the acceptance probability:
\begin{equation}
a(v,z) = \min\{1, \exp{[\Psi(z)-\Psi(v)]}\}, \label{e:acc}
\end{equation}
which is independent of discretization dimensionality up to numerical errors.

The pCN proposal in Eq. \eqref{e:pcn} can be improved by incorporating the data information in the proposal,
and following the idea of Langevin MCMC for the finite dimensional problems,  one can
derive the preconditioned Crank-Nicolson Langevin (pCNL) proposal:
\begin{equation}
(2+\delta)v=(2-\delta)z-2\delta \@C_0 \mathcal{D} \Psi(z) +\sqrt{8\delta }  w, \label{e:cnl}
\end{equation}
where $\delta\in[0,2]$, $w\sim N(\xi,\@C_0)$ and $\mathcal{D}$ is the gradient operator with respect to $z$. If we define $\rho(z,v)$ as following:
\begin{equation}\label{e:rho1}
\rho(z,v) = \Psi(z)+\frac{1}{2}\langle v-z,\mathcal{D}\Psi(z)\rangle +\frac{\delta}{4}\langle z+v,\mathcal{D}\Psi(z)\rangle+\frac{\delta}{4}\|    \@C_0^{1/2} \mathcal{D}\Psi(z) \|^2,
\end{equation}
then the acceptance probability is given by:
\begin{equation}
a(z,v) = \min\{1, \exp{(\rho(z,v)-\rho(v,z))  }  \}. \label{e:acc_pcnl}
\end{equation}

The pCNL algorithm is usually more efficient than the standard pCN algorithm as it takes advantage of the gradient information of the $\Psi(z)$.
However, the pCNL algorithm can not be used directly in our problem as $\Psi(z)$ includes the TV term which is not differentiable.
It is important to note that  $-2\delta \@C_0 \mathcal{D} \Psi(z) $ is the offset term only affecting the mean of the proposal, and we can replace $\mathcal{D} \Psi(z)$ with an alternative direction
$g(z)$, yielding proposal:
\begin{equation}
(2+\delta)v=(2-\delta)z-2\delta \@C_0 g +\sqrt{8\delta }  w, \label{e:pdprop}
\end{equation}
where $\delta\in[0,2]$ and $w\sim N(\xi,\@C_0)$.
Regarding the proposal given by Eq.~\eqref{e:pdprop}, we have the following theorem:
{\begin{theorem}\label{thm:db}
Assume that $\Psi$ satisfies Assumptions 6.1 in \cite{cotter2013mcmc}, and $g(z)$ is in the Cameron-Martin space associate with the Gaussian measure $\mu_0$.
Let $q(z,dv)$ be the conditional distribution defined by Eq.~\eqref{e:pdprop},
and define
 \[\eta(dz,dv)=q(z,dv)\mu^y(dz), \quad\eta^{\perp}(dz,dv)=q(v,dz)\mu(dv),\]
on on $X\times X$.
 We have that $\eta^{\perp}$ is equivalent to $\eta$ and
\begin{equation}
\frac{d\eta^{\perp}}{d\eta}(z,v)=\exp(\rho(z,v)-\rho(v,z)).
\end{equation}  
where
\begin{equation}
\rho(z,v) = \Psi(z)+\frac{1}{2}\langle v-z,g\rangle +\frac{\delta}{4}\langle z+v,g\rangle+\frac{\delta}{4}\|    \@C_0^{1/2} g \|^2.
\label{e:newrho}
\end{equation}
\end{theorem}
The proof of the theorem is provided in Appendix. It should be clear that Theorem~\ref{thm:db} 
implies that the MCMC algorithm with proposal~\eqref{e:pdprop} 
yields a well defined acceptance probability in the function space,
and as a result the chain satisfies the detailed balance condition in the function space and thus is ergodic. 
 Another  very important theoretical issue here is to estimate the spectral gaps and prove the geometric ergodicity of the algorithms in the infinite dimensional setting.
We note that there are some results on the spectral gaps of the standard pCN~\cite{hairer2014spectral} and the generalized pCN~\cite{rudolf2018generalization}.
It is an interesting problem to analyze if similar results can also be obtained for the present algorithm.}

Now we need to find a good direction $g(z)$.
In \cite{pereyra2016proximal} the authors use Moreau approximation to approximate the TV term in the Langevin MCMC algorithm.
Here we shall provide an alternative approach, determining the offset direction in the MCMC iteration using the primal dual algorithm.
The primal-dual algorithms are known
to be very effective in solving optimization problems involving TV regularization~\cite{chambolle2011first,chan1999nonlinear,chambolle2004algorithm},
and we hereby give a brief description of the primal dual method applied to our  problem.
Suppose that we want to solve
\begin{equation}
\min_{z\in X} \Psi(z)  = \Phi(z) +\lambda \|z\|_\mathrm{TV}.\label{e:opt1}
\end{equation}
Introducing a new variable $\bm\phi(x)=[\phi_1(x),\,\phi_2(x)]$ with
$\phi_1(\-x),\phi_2(\-x) \in L_2(\Omega)$ (we denote this as
$\bm\phi\in L_2^2(\Omega)$),
and we then rewrite the optimization problem~\eqref{e:opt1} as
\begin{align}\label{e:admmopt}
  \min_{z\in X,\phi\in L_2^2(\Omega)}  \quad   \Psi(z,\phi) &=\Phi(z)+\lambda \|\phi\|_{2,1} \\
  \textrm{s.t.}  \quad  \quad \nabla z &= \phi, \notag
\end{align}
where $\|\phi\|_{2,1} =\left(\|\phi_{1}(\-x)\|^2_{L_2(\Omega)}+\|\phi_2(\-x)\|^2_{L_2(\Omega)} \right)^{1/2}$.
The augmented Lagrangian for  Eq.\eqref{e:admmopt} is
\begin{equation}
\label{e:lagrange}
 \max_{\eta\in L_2^q(\Omega)} \min_{z\in X,\phi\in L_2^q(\Omega) } \quad L_\rho(z,\phi,\eta) = \Phi(z)+\lambda \|\phi\|_{2,1}
  + \<\eta, \nabla z-\phi \>
  +\frac{\rho}{2} \|\nabla z-\phi \|_2^2,
\end{equation}
where $\eta \in L_2^q(\Omega)$ is the dual variable or Lagrange multiplier, and $\rho>0$ is a constant called the penalty parameter.
The resulting dual problem is then solved
with the Alternating Direction Method of Multipliers (ADMM)~\cite{boyd2011distributed}:
\begin{subequations}
\begin{align}
z^{k+1}&=  \textrm{arg} \min_{z\in X}~~L_\rho(z,\phi^{k},\eta^{k}),  \label{e:admm1}\\
\phi^{k+1}&= \textrm{arg} \min_{\phi\in L^q_2(\Omega)}~~L_\rho(z^{k+1},\phi,\eta^{k}),   \label{e:admm2}\\
\eta^{k+1}&= \eta^{k}+ \rho(\phi^{k+1}-\nabla z^{k+1}).  \label{e:admm3}
\end{align}
\end{subequations}
The algorithm consists of a $z$-minimization step \eqref{e:admm1}, a $\phi$--minimization step \eqref{e:admm2}, and a dual ascent step  \eqref{e:admm3}.

 Our primal dual pCN(PD-pCN) algorithm is designed as follows. First we solve Eq.\eqref{e:lagrange} with the ADMM algorithm \ref{alg:mpcnl} obtaining the solution
  $(z^\star, \phi^\star, \eta^\star)$ and then we define
  \begin{equation}
\label{e:gz}
g(z) = \mathcal{T}_K\mathcal{D} L_{\rho} (z, \phi^\star, \eta^\star),
\end{equation}
where operator $\mathcal{T}_K$ is the projection of its input function onto the space spanned by the KL models $\{e_1,...c_K\}$ 
for a prescribed positive integer $K$. $K$ should be no greater  than the discretization dimensionality $N$. 
It should be clear that the function $g(z)$ computed with Eq.~\eqref{e:gz} is in the the Cameron-Martin space
of $\mu_0$.
The complete algorithm is given in Algorithm \ref{alg:mpcnl}. We note here that the proposed MCMC algorithm involves an optimization problem at the beginning  and the computational cost for solving this optimization is usually an order of magnitude lower that of the MCMC iterations.
{  A main limitation of this algorithm is that, when some hyper-parameters change, the optimization problem needs to be solved again, which makes it incompatible with  Metropolis within Gibbs~\cite{andrieu2003introduction} type of methods. 
It is  also worth noting here that, the main purpose of the proposed algorithm is to improve the sampling efficiency of the standard pCN algorithm while
maintaining its dimension independence property.  To this end a very interesting problem here is to incorporate the pCN framework with the aforementioned proximity based algorithm,
and compare the performance with the PD based one.}
\begin{algorithm}[!thb]
\caption{The Primal-dual pCN(PD-pCN) algorithm}
\label{alg:mpcnl}
\begin{algorithmic}[1]
\State Solve Eq.\eqref{e:lagrange} and denote the solution as $(z^\star, \phi^\star, \eta^\star)$
\State Let  $z^0 = z^\star$
   \For{$k=0,1,2,\cdots$}
   \medskip
    \State Propose v using Eq. \eqref{e:pdprop} and Eq. \eqref{e:gz}
    \State Draw $\theta \sim U[0,1]$;
    \State Compute $a(z,v)$ with Eq. \eqref{e:acc_pcnl} and Eq.~\eqref{e:newrho};
        \If {$\theta\leq a$ }
            \State $z^{k+1} = v$;
						\Else
						\State $z^{k+1}= z^{k}$;
        \EndIf
   \EndFor
\end{algorithmic}
\end{algorithm}
\section{Determining the hyperparameters}\label{sec:lambda}

Just like the deterministic inverse problems, it is an important issue to determine the TV regularization parameter $\lambda$ in the hybrid prior. 
In the Bayesian setting, the regularization parametter can be determined by the empirical Bayes (EB) approach~\cite{gelman2013bayesian}. 
Namely, the EB method seeks to maximize
\[\pi(\-y|\lambda) = \int \pi(\-y|u) \mu_\mathrm{pr}(du)
= \int \pi(\-y|u) \frac1{Z(\lambda)}\exp(-\lambda\|u\|_\mathrm{TV}) \mu_0(du)
\]
where $Z(\lambda)$ is the normalization constant.
{A difficulty here is that $Z(\lambda)$ is usually not known in advance and needs to be evaluated with another Monte Carlo integration. To address this issue, a stochastic proximal gradient method was proposed in \cite{pereyra2013computing,fernandez2018maximum}.
The method can efficiently estimate the regularization parameter $\lambda$ without the knowledge of $Z$. 
On the other hand, the method does require a suitable admissible set for $\lambda$.  
Here we provide a statistical approach to determine the admissible set of $\lambda$,
which is derived from the realized discrepancy method for model assessment proposed in~\cite{Gelman1996Posterior}. }
 
The basic idea of the method is to choose a function $D(\-y,\bm\theta)$ that measures the discrepancy between the measured data $\-y$ and the projected
observable $\bm\theta$, and for  the present problem we use the $\chi^2$ discrepancy, 
\begin{equation}
D(\-y,\bm\theta) =  \sum_{i=1}^d \frac{(y_i-\theta_i)^2}{\theta_i^2}. \label{e:chi2disc}
\end{equation}
Now knowing that $\theta = \~A u$, we can use this discrepancy to assess how well a specific choice of $u$ fits the data. 
The classical $p$-value based on the discrepancy $D(\-y,\bm\theta)$ is 
\begin{equation}
p_c(\-y,\bm\theta) =\P[ D(\-y,\bm\theta)>D(\tilde{\-y},\bm\theta)]
\end{equation}
where $\tilde{\-y}$ is the simulated data from model~\eqref{e:ytheta}.
In particular, for the discrepancy function given in Eq.~\eqref{e:chi2disc}, the $p$-value is simply, 
\begin{equation}
p_c(\-y,\bm\theta)=  1-F_{\chi_d^2}(D(\-y,\bm\theta)),
\end{equation}
where $F_{\chi_d^2}(\cdot)$ is the cumulative  distribution function of the $\chi^2$ distribution with the degree of freedom $d$. 
The classic $p$-value computed this way provides an assessment of how well a single estimate of $u$ fits the data $\-y$. 
The method can be extended to the Bayesian setting to assess the fitness of the posterior distribution to data. 
First recall that our prior distribution given by Eq.~\eqref{e:priorz} is specified by the parameter $\lambda$, 
 and as a result the posterior also depends on $\lambda$ and here we write the posterior as $\mu_\lambda^y(du)$
 to emphasize its dependence on $\lambda$. 
In the Bayesian setting, one can compute the posterior predictive $p$-value:  
\begin{equation}
p_b(\-y,\lambda) = \int  p_c(\-y,\bm\theta) p(\bm\theta|\-y,\lambda) d\bm\theta
= \int p_c(\-y,\~Au) \mu_\lambda^y(du),
\end{equation}
which is essentially the classical $p$-value averaged over the posterior distribution. 
 The posterior predictive $p$-value assesses the fitness of the posterior distribution to the data: 
intuitively speaking, larger value of $p_b$ indicates better fitness of the posterior to the data~$\-y$. 
However, one can not simply choose the value of $\lambda$ that yields the largest value of $p_b$, or, equivalently the best 
fitness to the data, as that may cause overfitting. 
In other words, if the posterior fits the data ``too well'', it often implies that the effect of the prior distribution is so weak 
that the posterior is dominated by the data. 
In the image reconstruction problem, this situation is greatly undesirable, as the problem is highly ill-posed 
and we need significant contribution from the prior distribution  to obtain good estimates of the unknown $u$. 
In this respect, we should choose $\lambda$ in a way that the effects of the prior and data are well balanced,
which should be indicated by an appropriate value of $p_b$. 
The suitable values of $p_b$ are certainly problem dependent, and in the present problem we suggest to choose $\lambda$ so that the resulting
value of $p_b$ is \emph{approximately} in the range of $0.1\sim0.7$. {Based on this, we choose the admissible set for $\lambda$ to be
\[
\Lambda = \{\lambda\geq0| 0.1\leq p_b(y,\lambda)\leq0.7\}.
\] 
The optimal value of $\lambda$ is then determined by using the method in \cite{fernandez2018maximum} within $\Lambda$. }
It is worth mentioning that methods using the data discrepancy to determine 
regularization parameters for Poisson data model have also been developed in the deterministic setting,
and we refer to \cite{bardsley2009regularization,bertero2010discrepancy} for further details.
 {It is important to note that the discrepancy principle may lead to over-smoothing in certain problems~\cite{hall1987common},
 which, however, may not cause issues in the proposed method as it just uses the discrepancy method 
to identify the admissible set of the regularization parameter while the actual value of it is determined with EB.} 
Finally we also note that, in addition to $\lambda$, the Gaussian distribution may also be subject to hyper-parameters, 
and in principle these hyper-parameters can be determined along with $\lambda$ using the proposed approach.  
However, here we choose not to do so for two reasons:
first determining multiple parameters may significantly increase the computational cost;
second, as the Gaussian distribution is merely used as a reference measure in our hybrid prior, the posterior distribution is not sensitive 
to it, and it usually suffices to choose these hyper-parameters based upon certain prior 
information (for example, historical data).

\section{Artifact detection using the posterior distribution}\label{sec:artifact}
In practical image reconstruction problems, due to the imperfection of methods or devices, a reconstructed image may {contain} what are not present
in the original imaged object. In this section we describe an application of the posterior distribution to detect artifacts in a reconstructed image.

Specifically, we consider the posterior distribution of the image at any given point $\-x$, which is denoted as $u_{\-x}$.
Consequently we can write the posterior distribution of $u_{\-x}$ as $\pi_{\-x}(u_{\-x}|\-y)$.
Next we consider the highest posterior density interval (HPDI) which is essentially the narrowest interval corresponding to a given confidence level.
More precisely, for an $\alpha\in[0,1]$, the $100(1-\alpha)\%$ HPDI
is defined as {\cite{pereyra2017maximum}}.
\[ C_\alpha=\{u(\-x) | \pi_{\-x}(u(\-x)|\-y)>\pi_\alpha\}
\]
where $\pi_\alpha$ is the largest constant satisfying
$\P[u_{\-x}| \pi_{\-x}(u_{\-x} |\-y)>\pi_\alpha]=1-\alpha$.
Now suppose that we have a reconstruct image $\hat{u}$ and we also write its value at $\-x$ as $\hat{u}_{\-x}$.
Next we shall estimate how large the credible level $(1-\alpha)$ must be so that the associated HDPI may contain $\hat{u}_{\-x}$.
That is, we compute the smallest value of $(1-\alpha)$ such that,
\[\hat{u}_{\-x} \in C_\alpha.\]
Intuitively speaking, the larger the computed credible level $(1-\alpha)$ is, the more likely the considered $\hat{u}_{\-x}$ is an artifact.
And we thus use the credible level $(1-\alpha)$ to measure how likely a point is an artifact, and
we can do this test for any point $\-x\in \Omega$.
Alternatively, the problem can also be formulated as a Bayesian hypothesis test with a fixed $\alpha$ (e.g. $\alpha=5\%$)~\cite{pereira2008can}: that is,
$\hat{u}(\-x)$ is regarded to be an artifact if it is not contained in the $(1-\alpha)$ HPDI for the prescribed value of $\alpha$.
However, it has been pointed out in  \cite{thulin2014decision} that performing hypothesis test with HPDI may cause certain theoretical issue and so here we choose not to use the  hypothesis test formulation.

It should be noted here that, in \cite{pereyra2017maximum,durmus2018efficient}, the authors utilize the highest posterior density (HDP) region to
test if a candidate image is likely to be a solution to the reconstruction problem.
The purpose of the present work differs from the aforementioned ones in that we want to identify
regions or pixels which are unlikely to be present in the original image, rather than to assess the entire image.

\section{Numerical results} \label{sec:results}

\begin{figure}[htb]
  \begin{subfigure}[t]{.32\textwidth}
    \centering
    \includegraphics[width=0.88\linewidth]{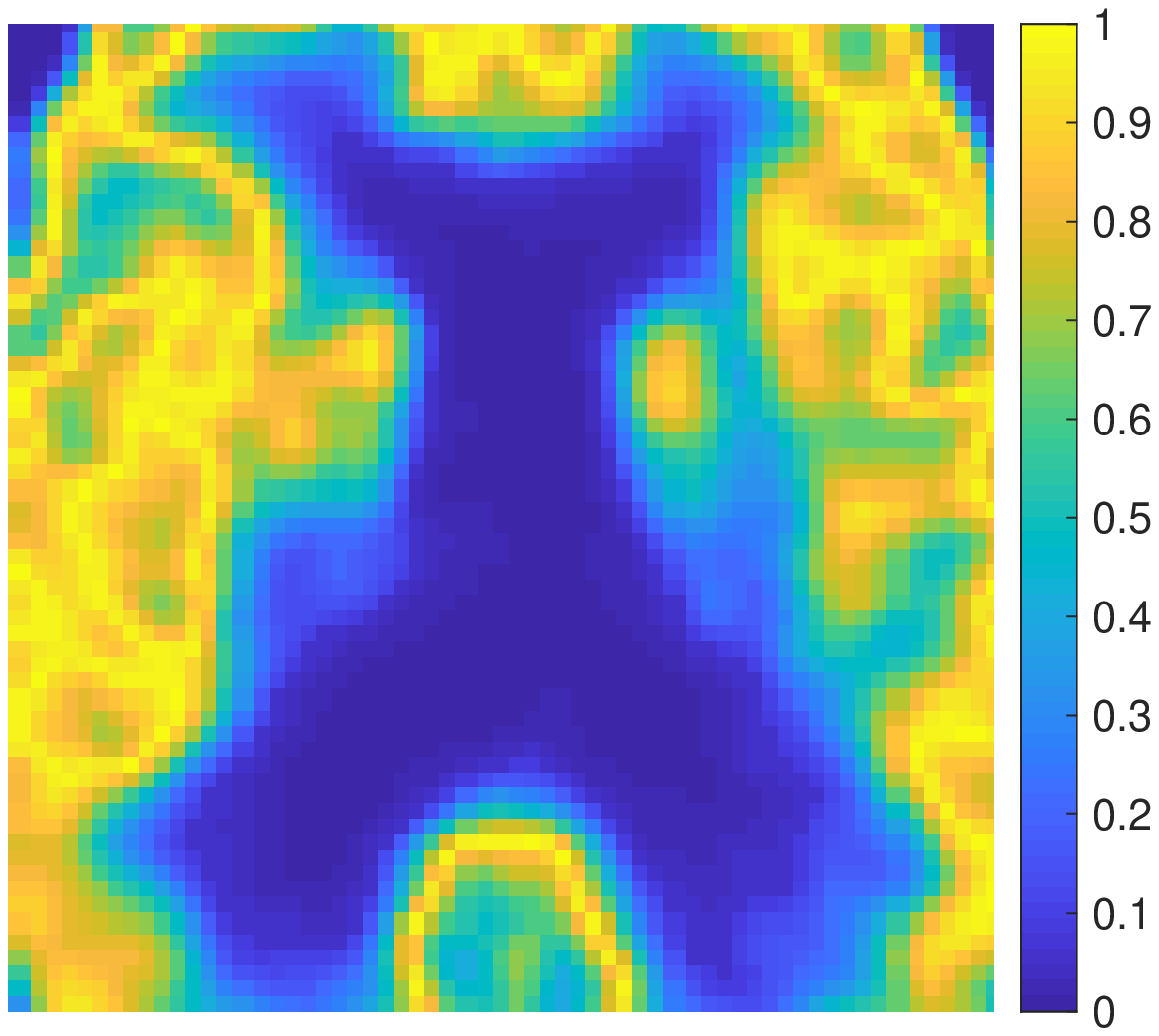}
  \end{subfigure}
  \begin{subfigure}[t]{.32\textwidth}
    \centering
    \includegraphics[width=0.88\linewidth]{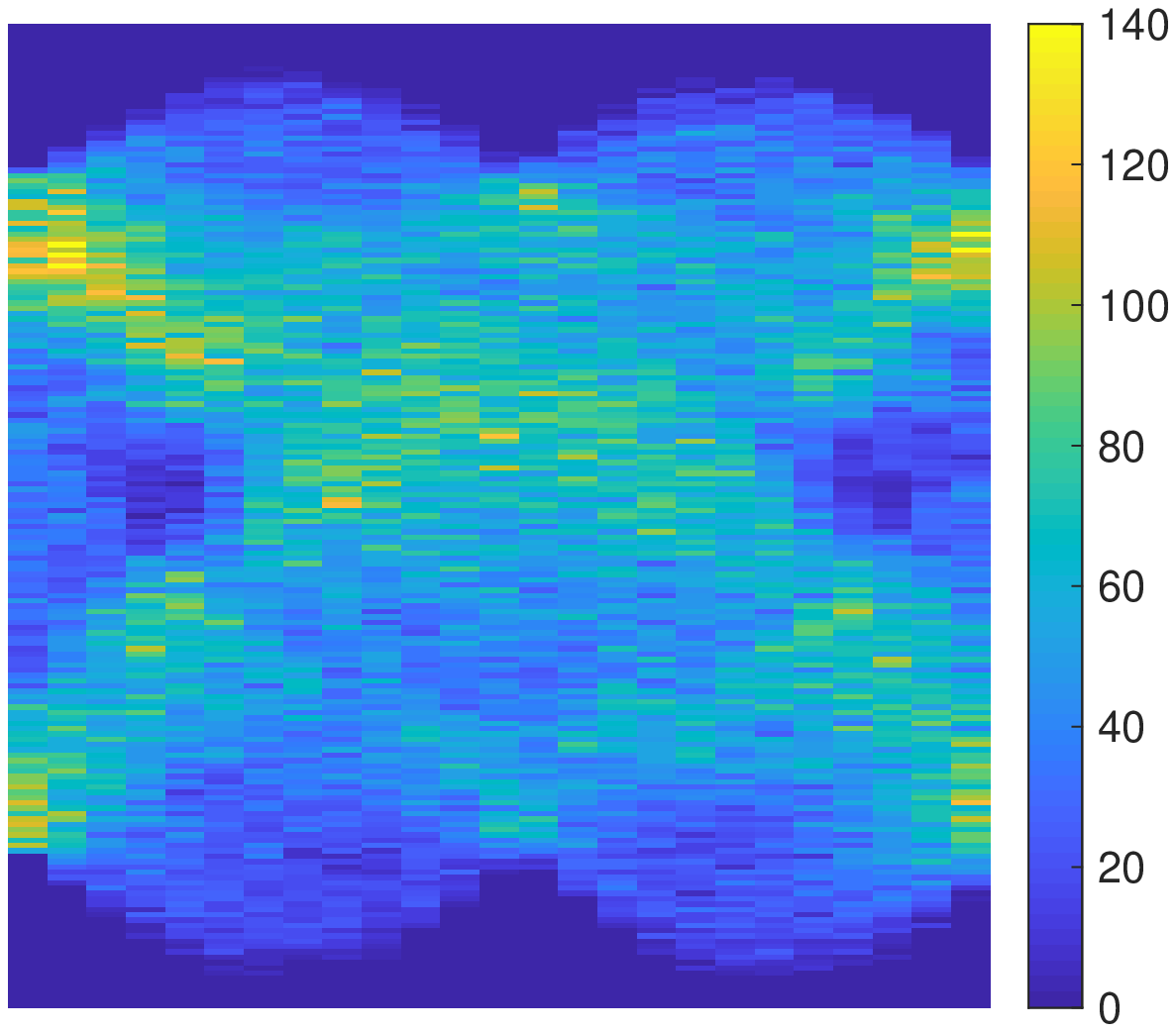}
  \end{subfigure}
  \begin{subfigure}[t]{.32\textwidth}
    \centering
    \includegraphics[width=0.88\linewidth]{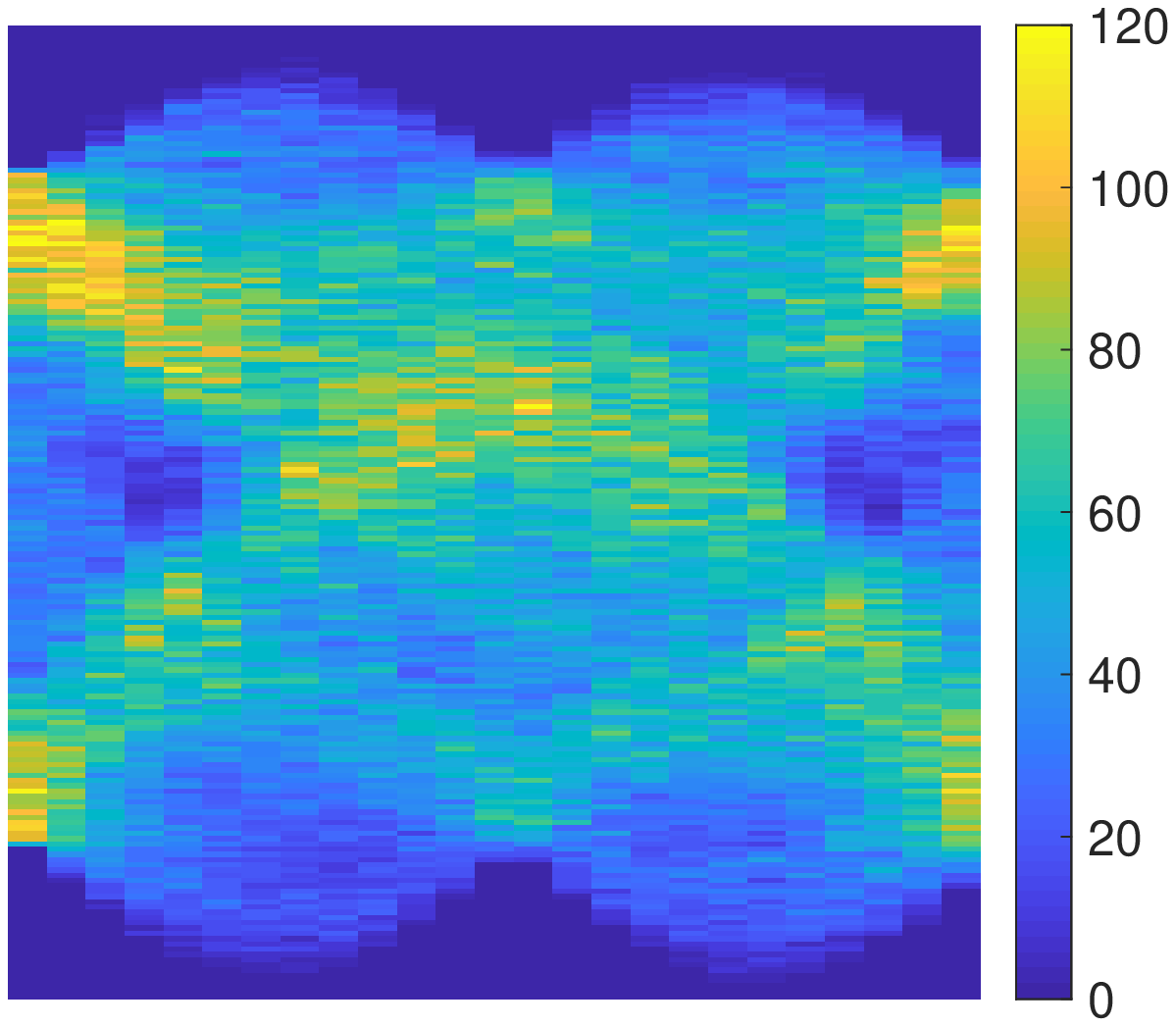}
  \end{subfigure}
  \caption{ Left: the true image. {Middle: the simulated data for $K=0.5$. Right: the simulated data for $K=1$}.}
  \label{f:truth-data}
\end{figure}

In this section we demonstrate the performance of the
proposed Bayesian framework, by applying it to a PET image reconstruction problem with synthetic data.
In particular the ground truth image ({Fig}.~\ref{f:truth-data}, left) is chosen from the Harvard whole brain atlas~\cite{petimage}.
 We let $\Omega=(0,1)^2$ and set the image size to be $128\times128$.
In the Radon transform we use 60 projections equilaterally sampled from 0 to $\pi$.
{In the numerical experiments, we consider two different noise levels: $K=0.5$ corresponding to a
higher noise level and $K=1$ corresponding to a lower noise level.
 The test data, shown in  {Figs}.~\ref{f:truth-data}, are randomly simulated by plugging the true image
into the Radon transform and the Poisson distribution~\eqref{e:ytheta} where the two aforementioned noise level $K$. }
In the Bayesian inference, we use the hybrid prior distribution where the Gaussian part is taken to be zero mean and covariance:
\begin{equation}\label{e:cov}
  K(\-x,\-x') = \gamma \exp \left[  -\frac{\|\-x-\-x'\|_1}{d} ,   \right]
\end{equation}
{where  $d$ is taken to be $10^{-3}$ and $\gamma$ is $2$.
The regularization parameter $\lambda$
are determined by using the method
presented in Section~\ref{sec:lambda},
and details will be discussed in next section}.

\subsection{Determining parameter $\lambda$}
As is discussed at the beginning of the section,  the prior parameter $\lambda$ is determined with the realized discrepancy method discussed in
Section~\ref{sec:lambda}.  We here provide some details on the issue.
Specifically we test five different values of $\lambda$ { for $K=0.5$:  $\lambda=0,1,2,3,4,5$, and using the method discussed in Section~\ref{sec:lambda} we
compute the corresponding posterior predictive $p$-value for each value of $\lambda$, shown in Table~\ref{t:pvalue2}. 
Similarly we also test 5 values of $\lambda$ for $K=1$ and the results are shown in 
Table~\ref{t:pvalue}. }
We can see from the table that, as $\lambda$ increases, the resulting $p$-value decays. 
These results agree
well with our expectation that as $\lambda$ becomes larger, the prior distribution becomes stronger, and as a result the $p$-value
which assesses the fitness of the posterior to the data becomes smaller.
We also compute the PSNR of the posterior distribution computed with all these $\lambda$ values, and the results are also given in the table.
We can see here that, for both very large and very small $p$-values, the associated posterior means are of rather poor quality in terms of  PSNR.
That is, when the $p$-value is too large,
the posterior distribution overfits the data, and when it is too small,
the posterior underfits the data;
both cases lead to a poor performance of the inference, and so we must choose a proper $p$-value that represents
a good balance of the prior and the data.
{Bases on the test results, for $K=0.5$ we choose $\Lambda=[1,\,3]$ and for $K=1$ we choose $\Lambda=[0.5,\,2]$.
By optimizing $\lambda$ within the identified intervals we obtain
$\lambda=2.4$ for $K=0.5$ and $\lambda=1.2$ for $K=1$.}

\begin{table}[htb]
  \centering
     \begin{tabular}{|c|c|c|c|c|c|c|}
       \hline
        $\lambda$ & 0 & 1 & 2 & 3 & 4 & 5 \\
        \hline
        $\textrm{$p_b$}$ & 0.99 & 0.74 & 0.32 &0.08 &0.0074&0.0004\\
        \hline
       $\textrm{PSNR}$ & 15.69 & 18.85 & 20.31 & 20.04 & 18.79 & 18.58 \\
       \hline
     \end{tabular}
  \caption{{(K=0.5) The posterior predictive $p$-value ($p_b$) and the PSNR of the resulting posterior mean for different values of $\lambda$}.}
  \label{t:pvalue2}
\end{table}

\begin{table}[htb]
  \centering
     \begin{tabular}{|c|c|c|c|c|c|c|}
       \hline
        $\lambda$ & 0 & 0.5 & 1 & 2 & 3 & 4 \\
        \hline
        $\textrm{$p_b$}$ & 0.99 & 0.82 & 0.28 &0.04 &0.0034&0.0005\\
        \hline
       $\textrm{PSNR}$ & 18.21 & 21.90 & 21.99 & 20.66 & 19.76 & 19.27 \\
       \hline
     \end{tabular}
  \caption{(K=1) The posterior predictive $p$-value ($p_b$) and the PSNR of the resulting posterior mean for different values of $\lambda$.}
  \label{t:pvalue}
\end{table}

%


\subsection{Convergence with respect to discretization dimensionality}
\begin{figure}[!htb]
    \centering
   \centerline{\includegraphics[width=0.5\linewidth]{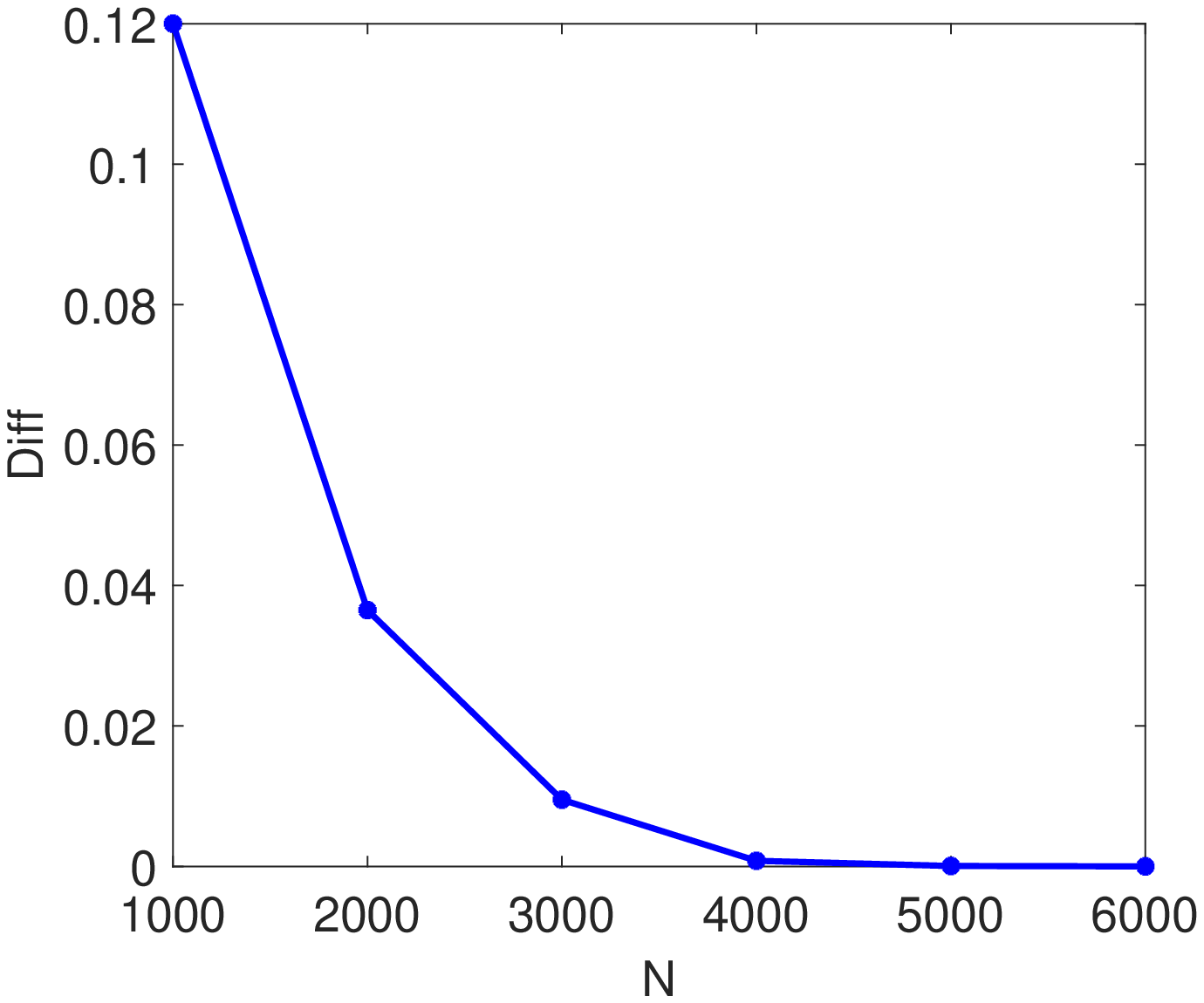}
     \includegraphics[width=0.5\linewidth]{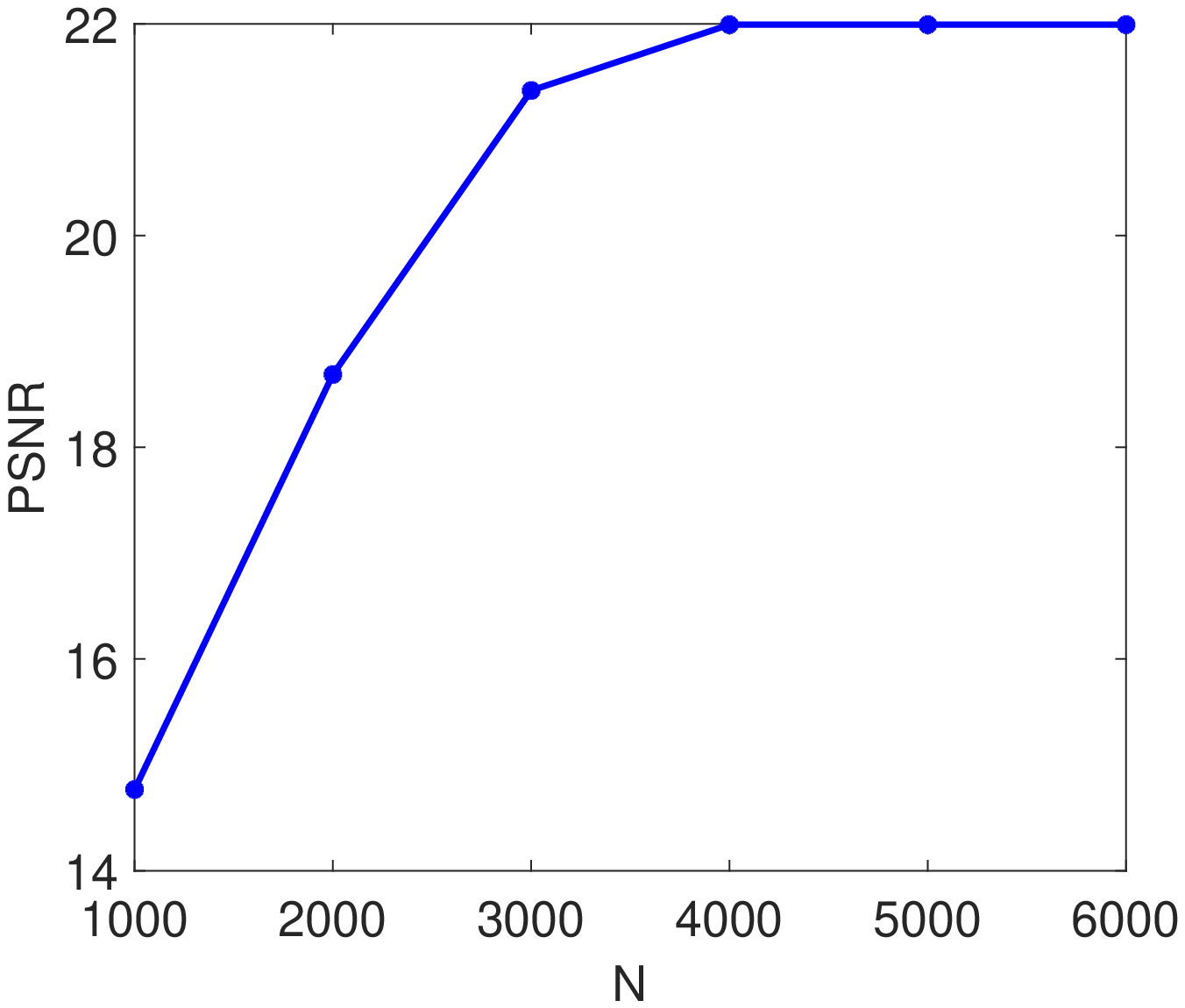}}
    \caption{ Left: the convergence of the posterior mean. Right: the PSNR of the posterior mean as a function of $N$.}
  \label{f:N}
\end{figure}

\begin{figure}[!htb]
    \centering
   \centerline{\includegraphics[width=0.8\linewidth]{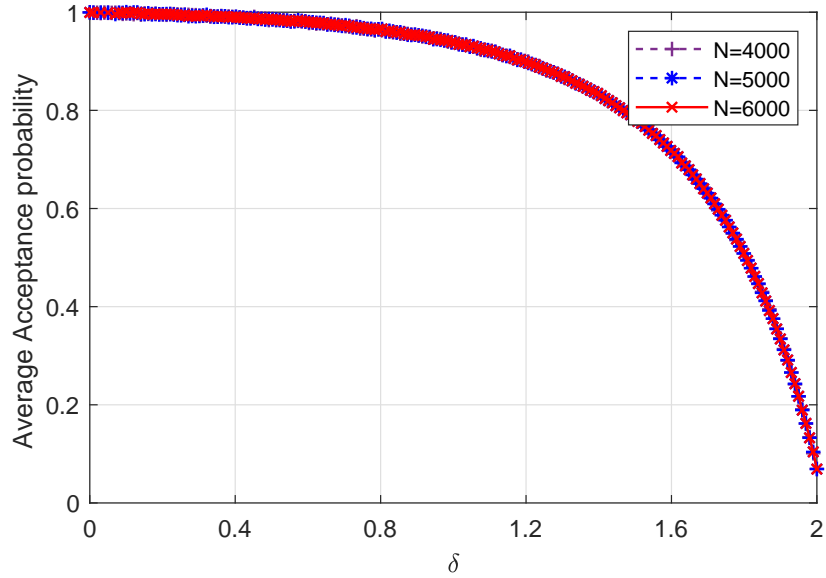}}
    \caption{ The average acceptance probability plotted against $\delta$ for $N=4000,\,5000,\,6000$.}
  \label{f:acc}
\end{figure}

In the numerical implementation we represent the unknown $z$ using the truncated KL expansion with $N$ KL-modes.
First we shall demonstrate that the posterior distributions converges with respect to the discretization dimensionality $N$.
We here use the case $K=1$ as an example. 
We perform the proposed PD-pCN MCMC simulation and compute the posterior means with six different values of $N$: $N_i=i\times10^3$ for $i=1...6$.
We note here that, in all the MCMC simulations performed in this section, we fix the number of samples to
be $5\times 10^5$ with additional $0.5\times10^5$ samples used in the burn-in step,
and also, in all the simulations the stepsize $\beta$ has been chosen in a way that the resulting acceptance probability is around $25\%$.
We then compute the $L_2$ norm of the difference between the posterior mean with $N=N_i$ and that with $N=N_{i+1}$ for each $i=1...6$:
\begin{equation}
\mathrm{Diff}= \int_\Omega (\hat{u}_{N_i}(\-x)-\hat{u}_{N_{i+1}}(\-x))^2 d\-x
\end{equation}
where $\hat{u}_{N_i}$ is the posterior mean of $u$ computed with $N_i$ KL modes.
We plot the $L_2$ difference against the discretization dimensionality $N$ in Fig.~\ref{f:N}.
One can see from the figure that the difference decrease as $N$ increase and
the difference becomes approximately zero for $N=5000$ and $N=6000$,
indicating the convergence of the posterior mean with respect to $N$.
For each posterior mean $\hat{u}_{N_i}$, we also compute its peak signal-to-noise ratio (PSNR)~\cite{hore2010image}, a commonly used metric
of the quality of a reconstructed image.
We show the PSNR results in Fig.~\ref{f:N} ({right}), and
the figure shows that the PSNR increases as $N$ increases from 1000 to 4000, and remains approximately constant from 4000 to 6000,
suggesting that increasing the discretization dimensionality can improve the inference accuracy until the posterior converges,
and so it is important to use sufficiently large discretization dimensionality in such problems.
{Next to further demonstrate that the proposed PD-pCN MCMC algorithm is independent of discretization dimensionality, 
we perform the MCMC simulation with different values of $\delta$ which is the parameter controlling the step size of the algorithm. In Fig~\ref{f:acc} we plot the average acceptance probability as a function of $\delta$ for 
three different values of $N$, and one can see that the acceptance probabilities 
under different discretization dimensionality agree
well with either other, indicating that the acceptance probability of the algorithm is independent of 
the discretization dimensionality $N$. }
In the rest of the work, we fix $N=6000$.

\begin{figure}
    \includegraphics[width=0.49\linewidth]{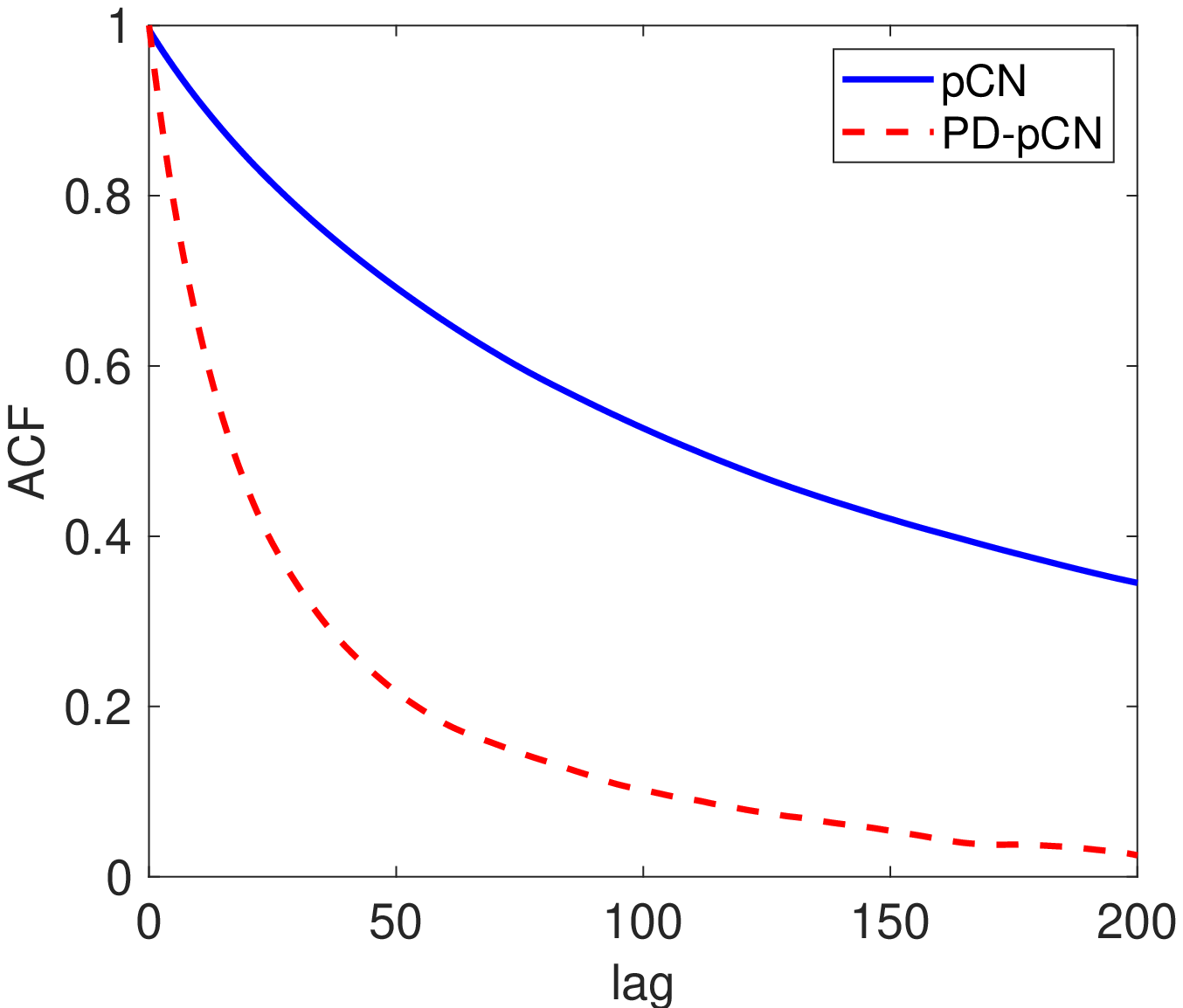}
    \includegraphics[width=0.49\linewidth]{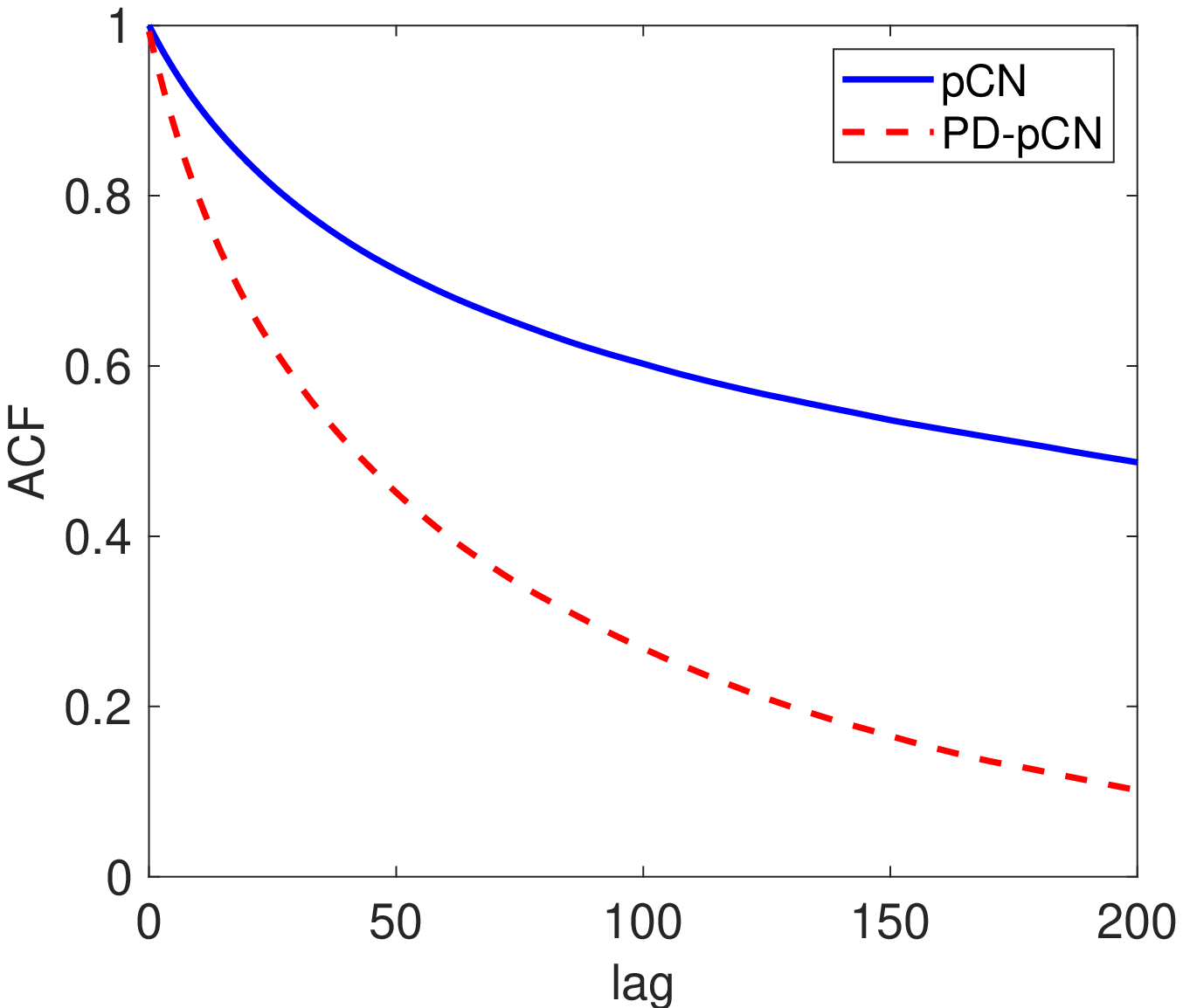}



  \caption{{(K=0.5)} The ACF of {the fastest~(left) and the slowest~(right) components of the samples} drawn by the pCN and the {PD-pCN} methods.}
  \label{f:acfKsmall}
\end{figure}

\begin{figure}
    \includegraphics[width=0.49\linewidth]{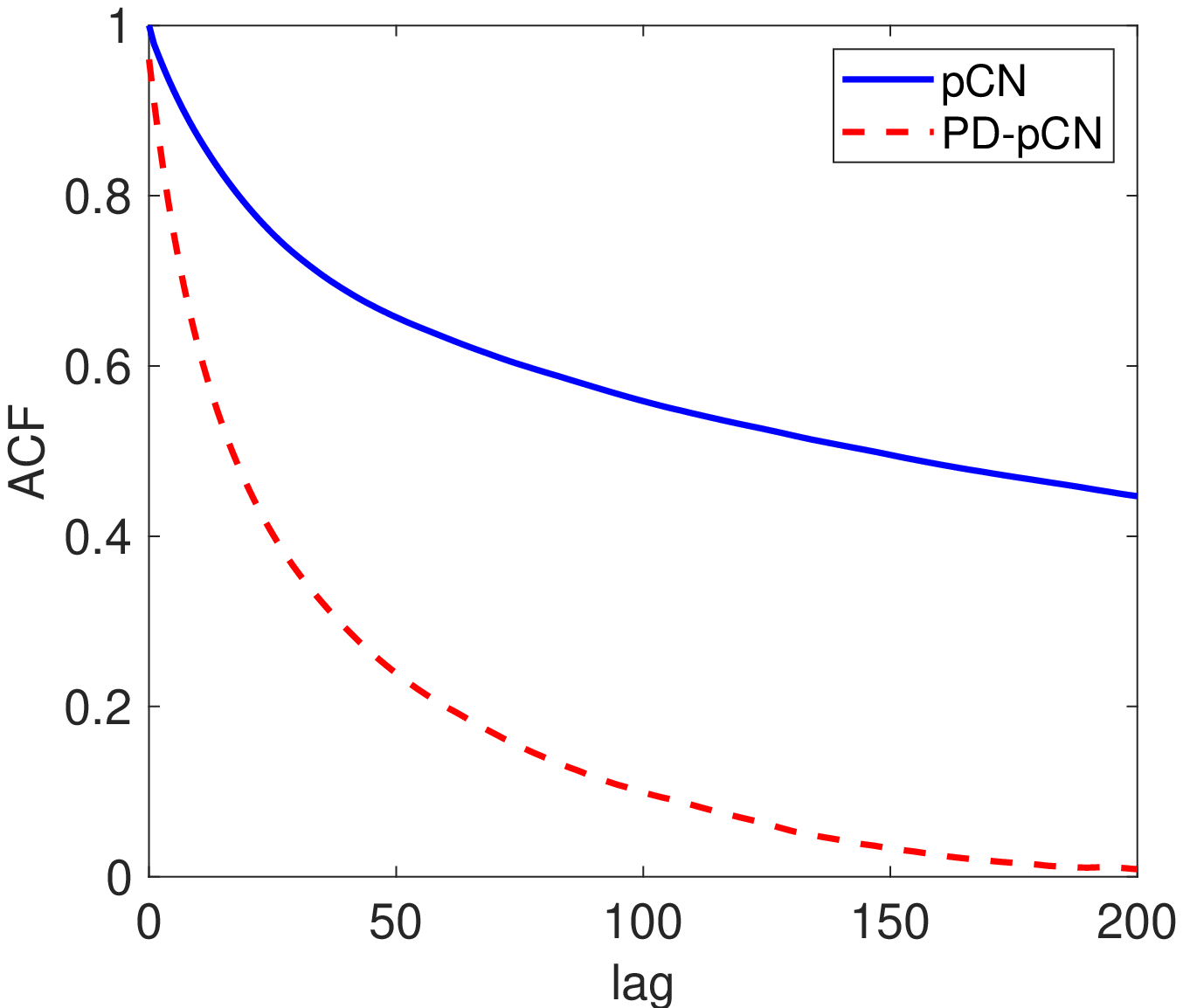}
    \includegraphics[width=0.49\linewidth]{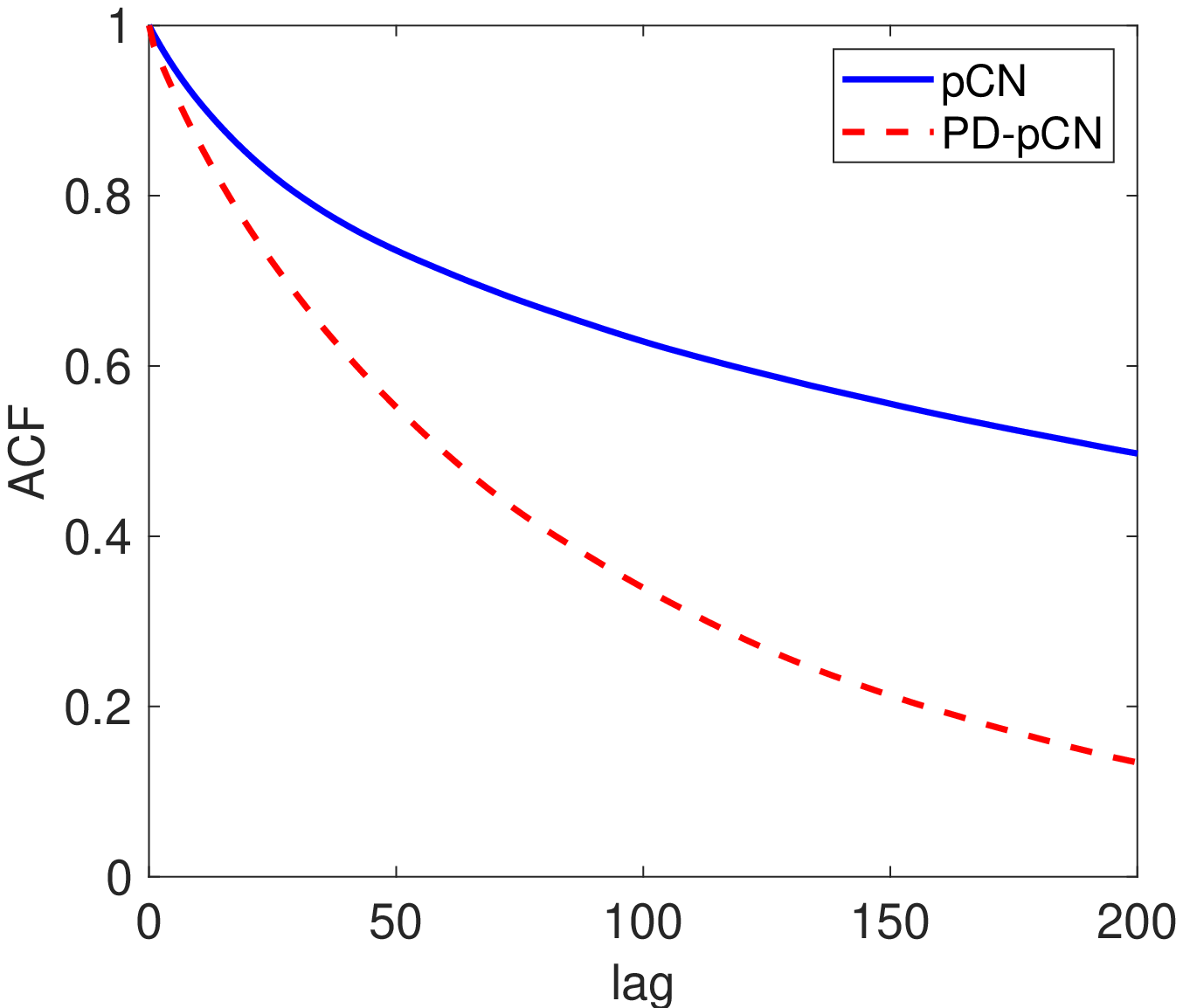}



  \caption{{(K=1)} The ACF of {the fastest~(left) and the slowest~(right) components of the samples} drawn by the pCN and the {PD-pCN} methods.}
  \label{f:acfKbig}
\end{figure}

\begin{figure}
    \centering
    \includegraphics[width=0.99\linewidth]{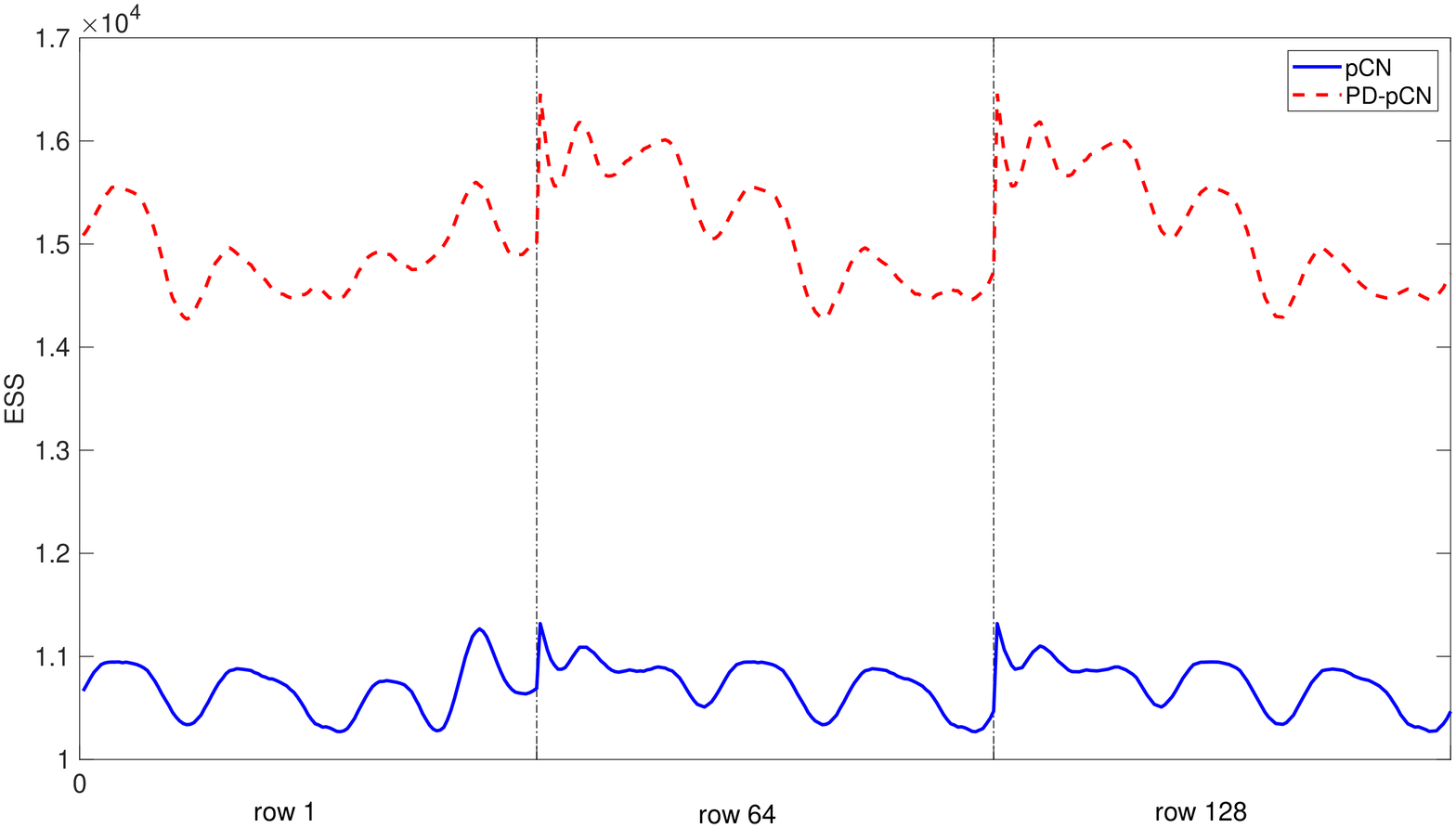}
    \caption{  ESS of {the} pCN and {the PD-pCN} algorithms}
  \label{f:ess}
\end{figure}

\subsection{Sampling efficiency of the PD-pCN algorithm}
Next, we shall compare the performance of the proposed PD-pCN algorithm and the standard pCN.
We draw $5\times10^5$ samples from the posterior distribution using both the standard pCN and the proposed PD-pCN algorithms.
We reinstate that in both algorithms we have chosen the step size so that the resulting acceptance probability is around $25\%$.
{In particular, to achieve the sought acceptance probability,  the values of the stepsize parameter $\beta$ in pCN are taken to be $0.04$ (for $K=1$) and $0.09$ (for $K=0.5$);
the values of the stepsize parameter $\delta$ in PD-pCN are $0.18$ (for $K=1$) and $0.23$ (for $K=0.5$).}
 {The total computational time is around 12 hours in a workstation with a 6-core 2.50 GHZ processor.}
We compute the auto-correlation function (ACF) of the samples generated by the two methods at all the grid points,
and {we show the ACF at the points with the fastest and the slowest convergence rates, 
in Fig.~\ref{f:acfKsmall} (K=0.5) and Fig.~\ref{f:acfKbig} (K=1).}
One can see from the figures that at both points the ACF of the propose PD-pCN method decays much faster than that of the standard pCN. {To further compare the performance,  we compute the effective sample size(ESS) which is defined as,
\[
\textrm{ESS}=\frac{N}{1+2\tau} ,
\]
where $\tau$ is the integrated auto-correlation time and $N$  is total sample size.}
In Fig.~\ref{f:ess},  we compare the ESS at three chosen rows from left to right in the image, namely row 1,  64 and 128,
for $K=1$. 
Just as the ACF, the results show that the PD-pCN algorithm achieves much higher ESS than the standard pCN.
We have also examined the ESS for $K=0.5$, where the results are qualitatively similar to those three shown in Fig.~\ref{f:ess},
and so we omit those results.

\subsection{The inference results}
\begin{figure}
  \begin{subfigure}[t]{0.24\textwidth}
    \centering
    \includegraphics[width=0.9\linewidth]{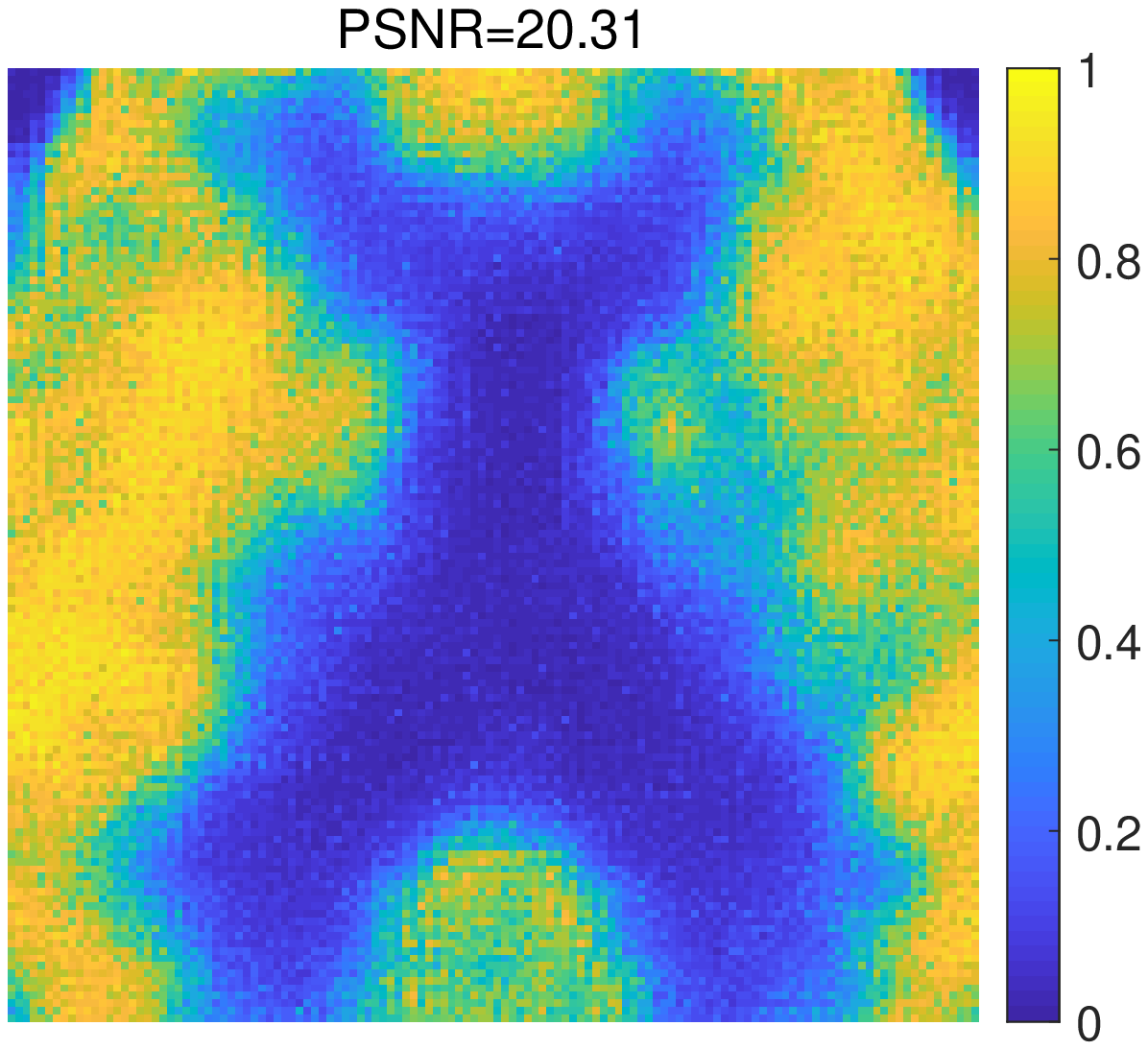}\label{f:k5tv2}
  \end{subfigure}
  \begin{subfigure}[t]{0.24\textwidth}
    \centering
    \includegraphics[width=0.9\linewidth]{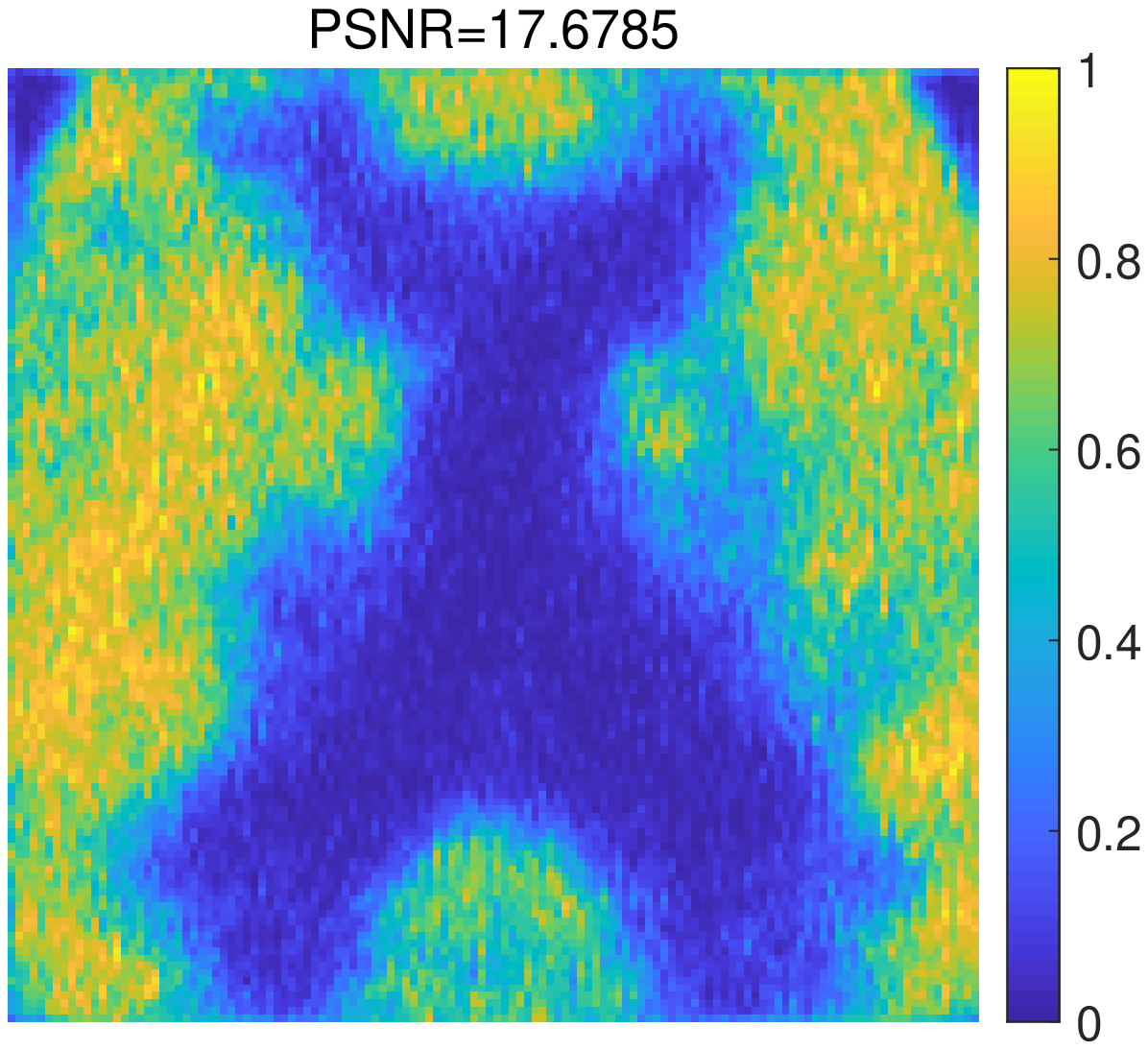}\label{f:k5tv0}
  \end{subfigure}
  \begin{subfigure}[t]{0.24\textwidth}
    \centering
    \includegraphics[width=0.9\linewidth]{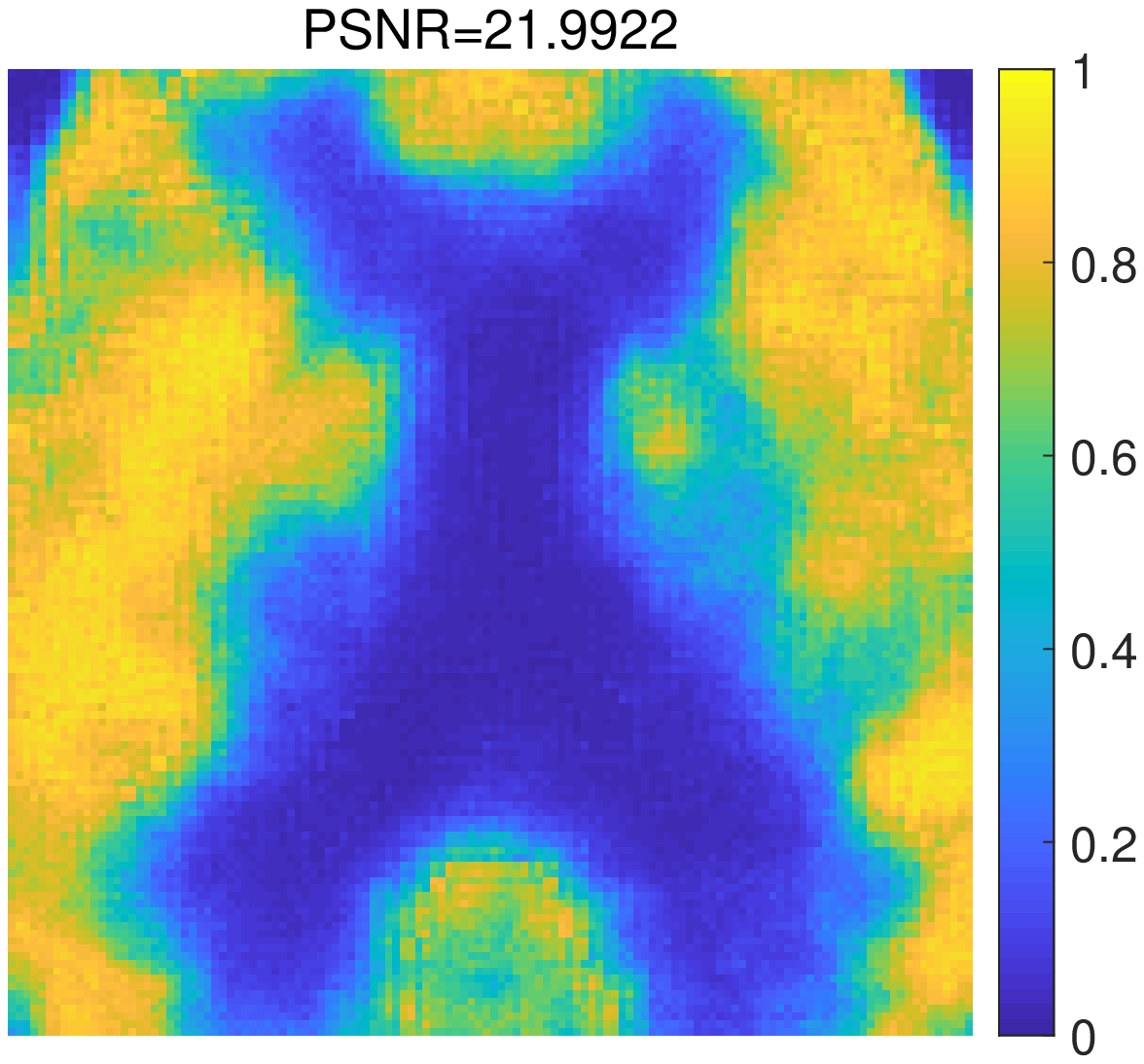}\label{f:k1tv1}
  \end{subfigure}
  \begin{subfigure}[t]{0.24\textwidth}
    \centering
    \includegraphics[width=0.9\linewidth]{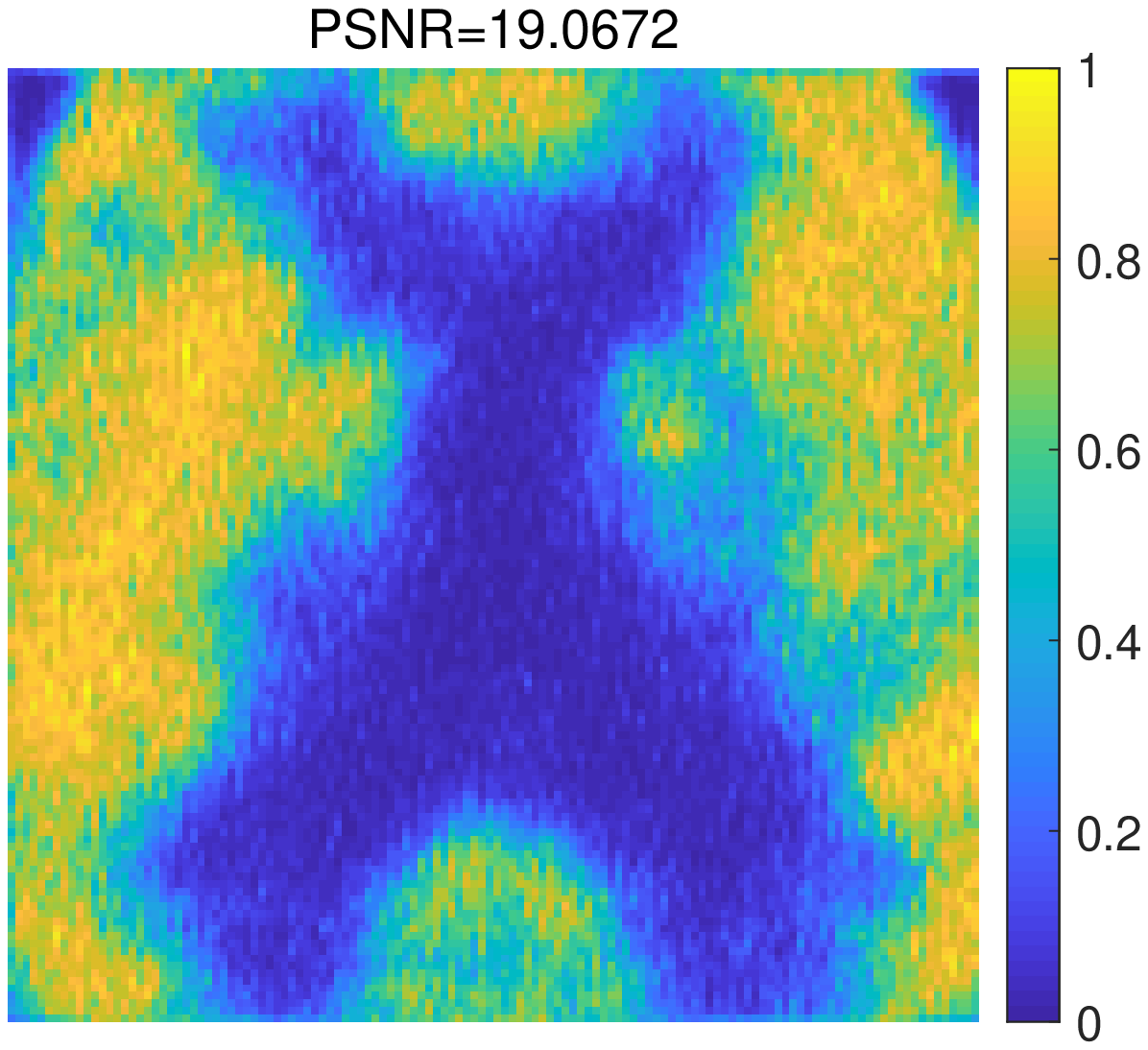}\label{f:k1tv0}
  \end{subfigure}

  \begin{subfigure}[t]{0.24\textwidth}
    \centering
    \includegraphics[width=0.9\linewidth]{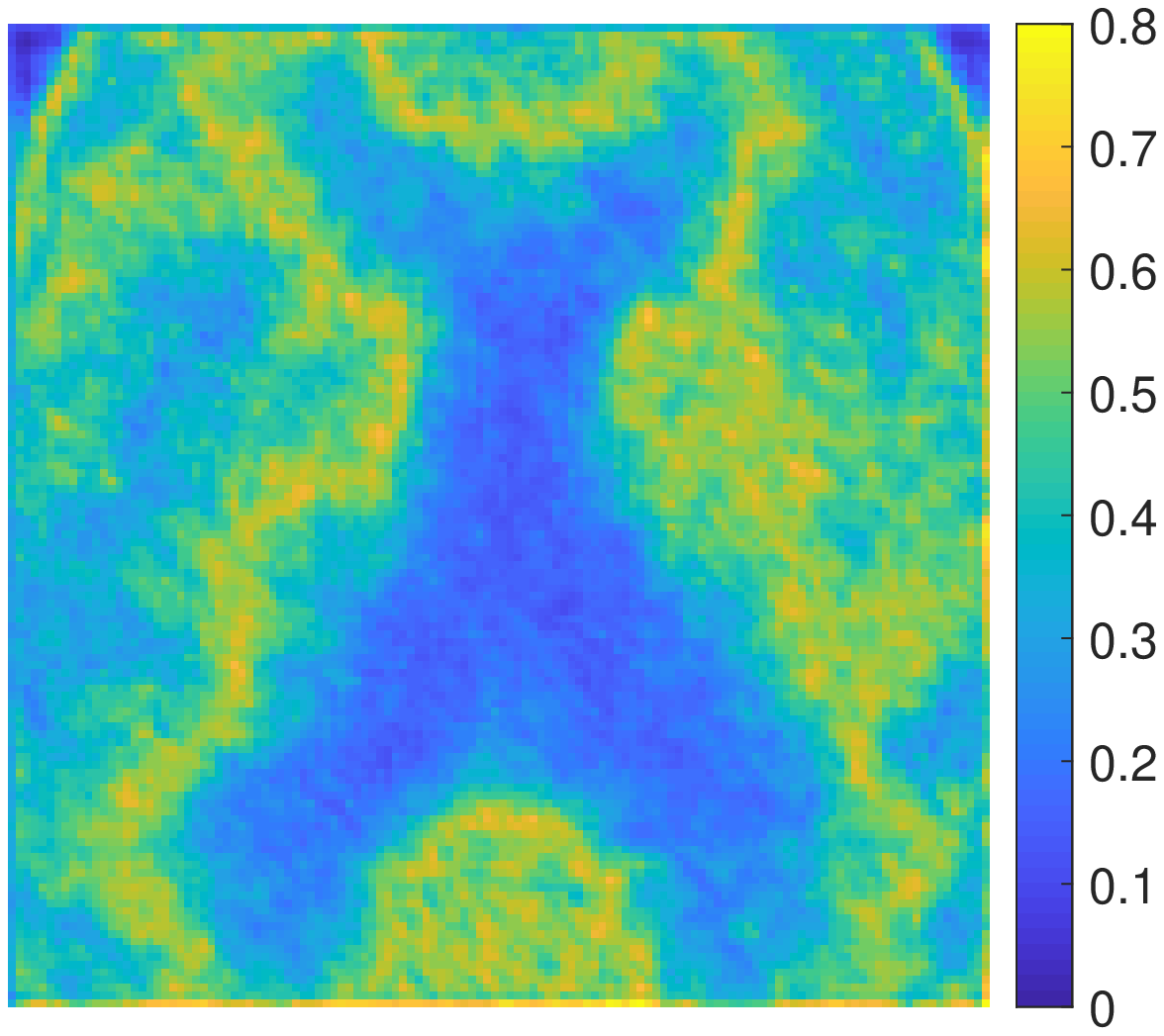}
    \caption{$K=0.5$, TG}
    \label{f:k5tv2ci}
  \end{subfigure}
  \begin{subfigure}[t]{0.24\textwidth}
    \centering
    \includegraphics[width=0.9\linewidth]{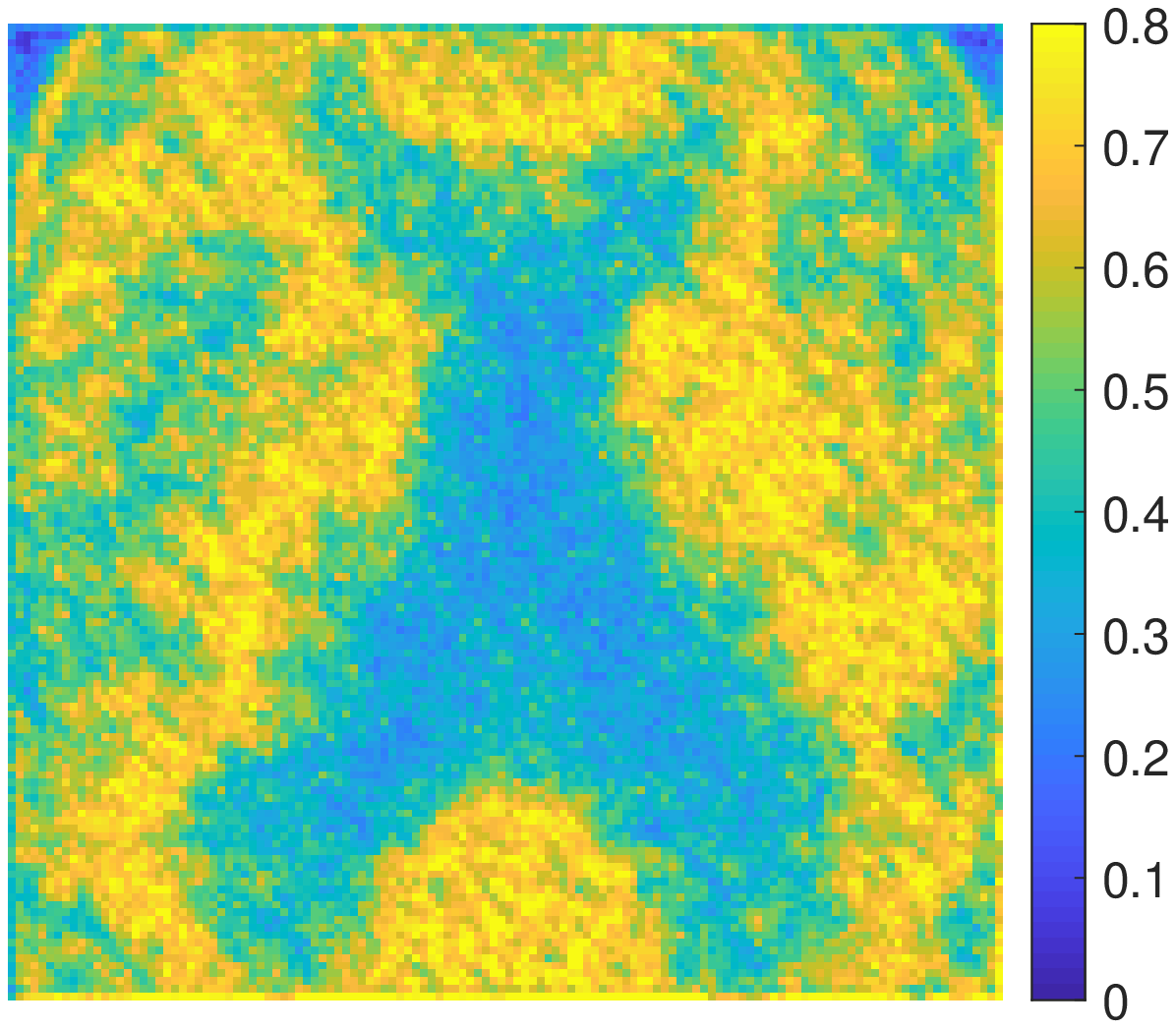}
    \caption{$K=0.5$, Gaussian}
    \label{f:k5tv0ci}
  \end{subfigure}
  \begin{subfigure}[t]{0.24\textwidth}
    \centering
    \includegraphics[width=0.9\linewidth]{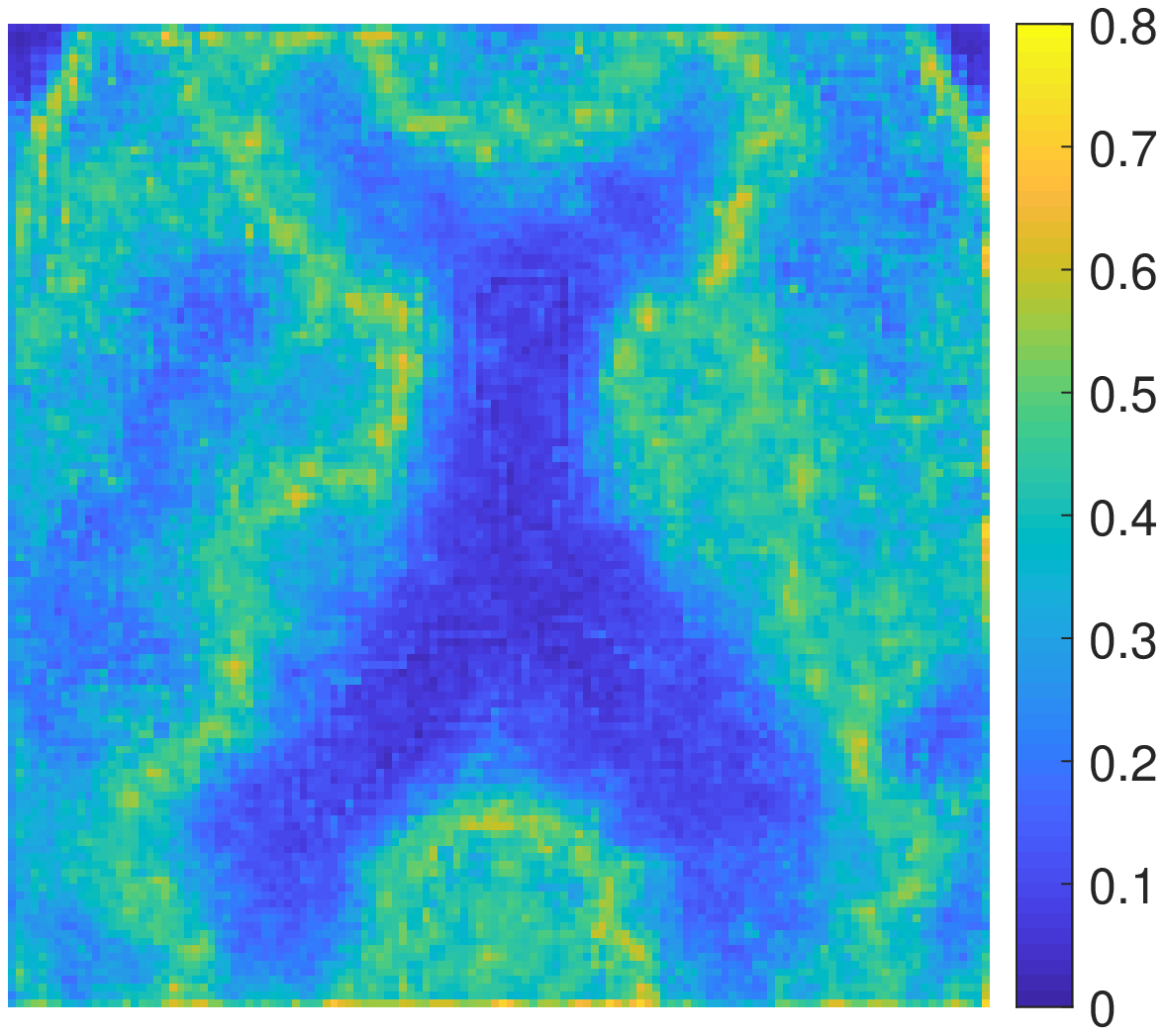}
    \caption{$K=1$, TG}
    \label{f:k1tv1ci}
  \end{subfigure}
  \begin{subfigure}[t]{0.24\textwidth}
    \centering
    \includegraphics[width=0.9\linewidth]{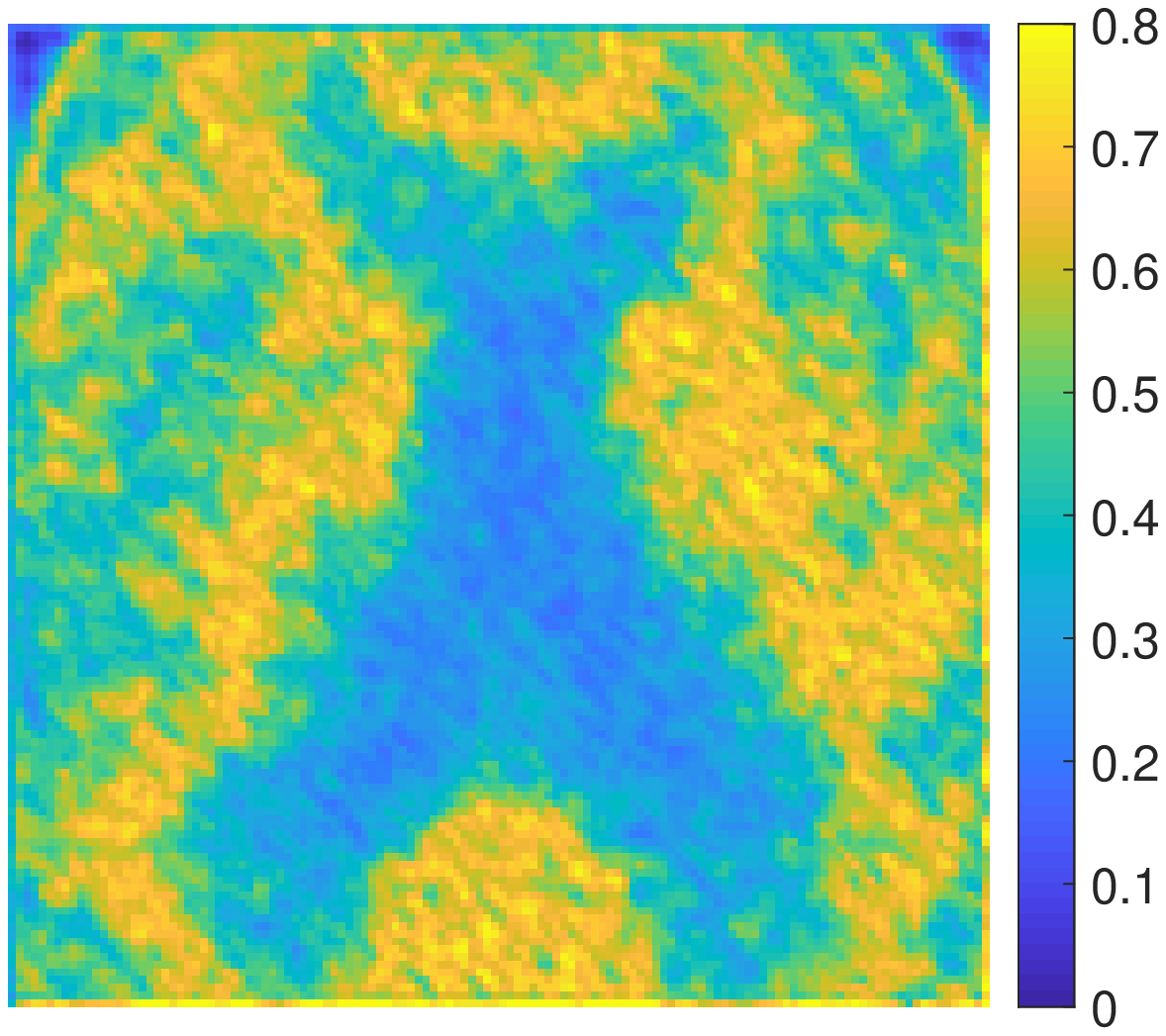}
       \caption{$K=1$, Gaussian}
       \label{f:k1tv0ci}
  \end{subfigure}
  \caption{{The posterior results of the TG prior and the Gaussian prior. The top figures show the posterior mean and the bottom ones shows the corresponding interval width of 95\% HPDI.}}
  \label{f:scale1}
\end{figure}

To illustrate the inference results, we compute the posterior mean of the TG prior, which is regarded as a point estimate of the image.
As is mentioned earlier, a main advantage of the Bayesian method is its ability to quantify the uncertainty in the reconstruction and to this end,
 we use the width of the (pointwise) 95\% HPDI as a metric of the posterior uncertainty~(intuitively speaking the wider the HPDI is, the more uncertainty there is). We plot these posterior results in {Figs}.~\ref{f:scale1}: { the posterior mean and the 95\% HPDI for $K=0.5$ and $K=1$ are shown in Fig.~\ref{f:k5tv2ci} and Fig.~\ref{f:k1tv1ci} respectively.}
 As a comparison, we also compute the posterior mean as well as the 95\% HPDI width, for the Gaussian prior corresponding to setting $\lambda=0$ in the TG prior, and the results are also shown in  {Fig.~\ref{f:k5tv0ci} and Fig.~\ref{f:k1tv0ci}} .
 The figures show that the posterior mean obtained with the TG prior is clearly of better quality than that of the Gaussian prior,
 suggesting that including the edge-preserving TV term significantly  improves the performance of the prior.
 It is worth noting here that, the Gaussian prior used here is not optimized for the best performance, and the performance can be potentially improved
 by using some carefully designed Gaussian priors, for example, \cite{calvetti2007gaussian}.

\begin{figure}[!htp]
  \begin{subfigure}[t]{.5\textwidth}
        \centering
        \includegraphics[width=.92\linewidth]{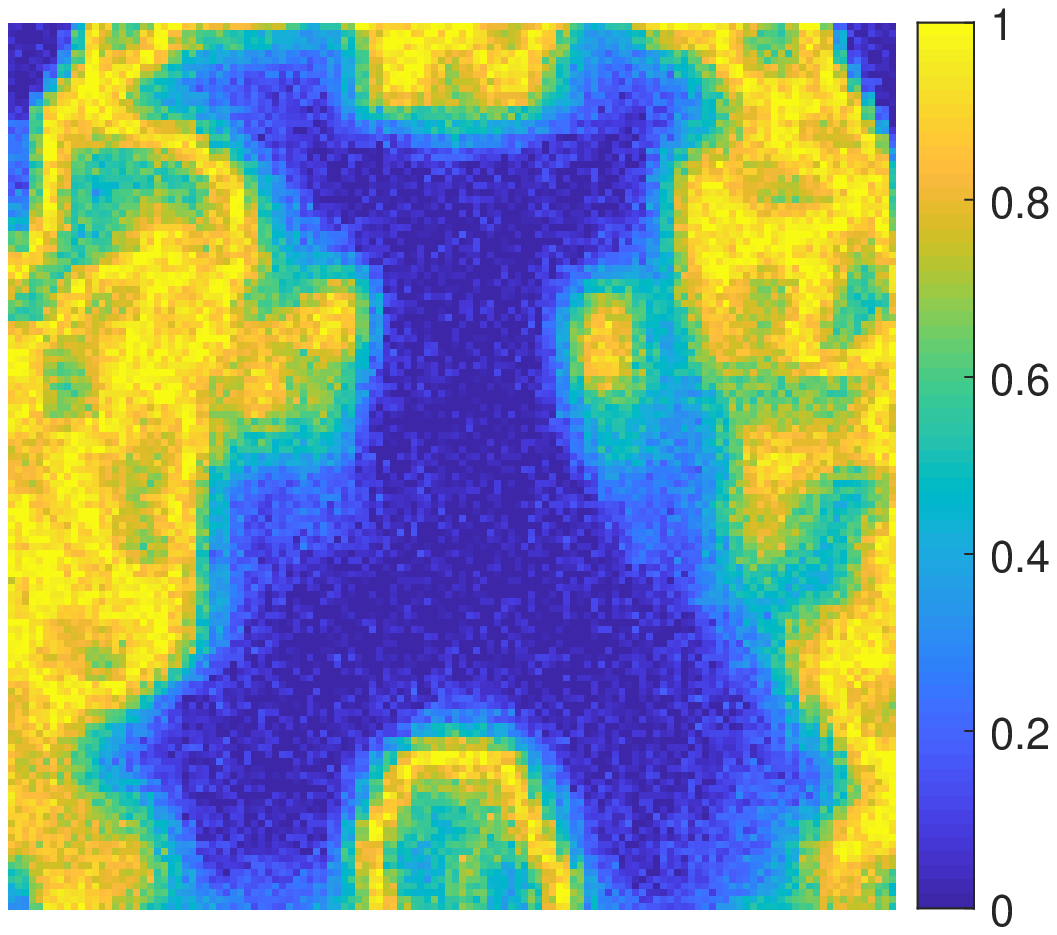}
\caption{Surrogate test image one.}
        \label{f:hypo0}
  \end{subfigure}
	~~~~
  \begin{subfigure}[t]{.5\textwidth}
        \centering
        \includegraphics[width=.92\linewidth]{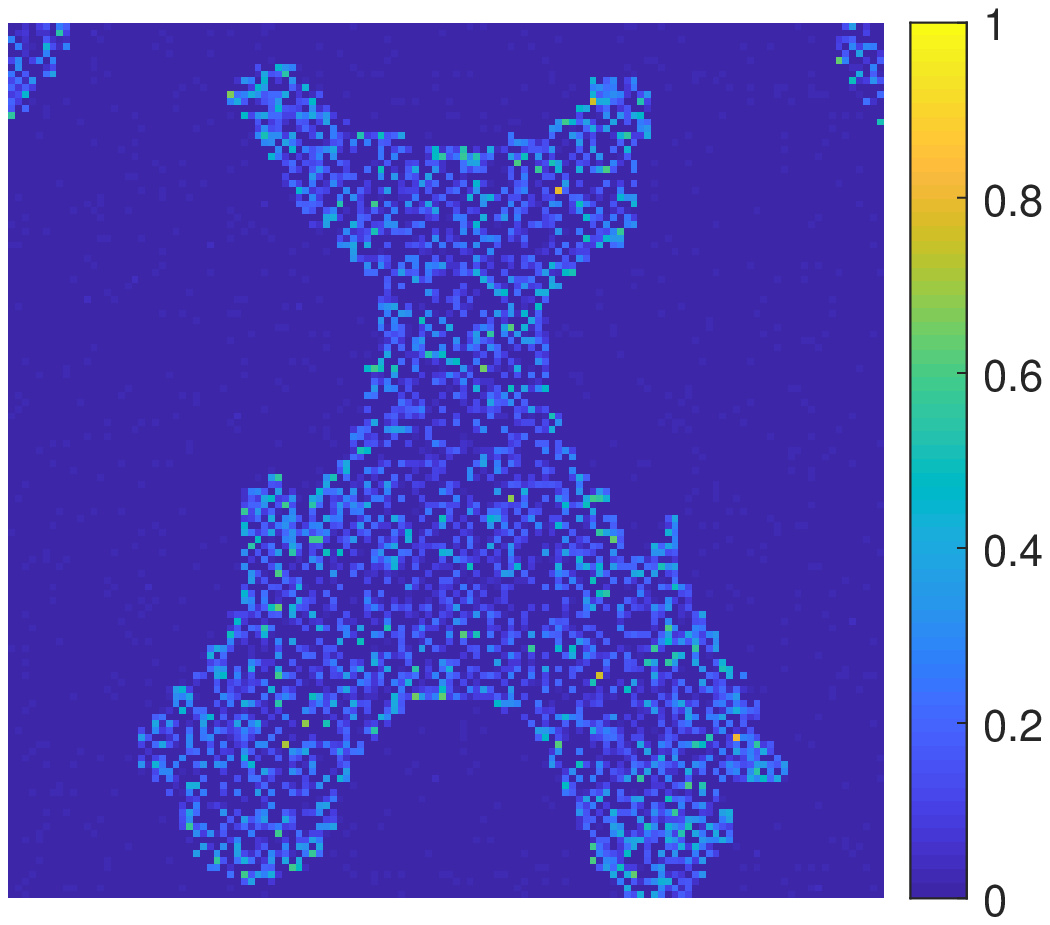}
\caption{The credible level $(1-\alpha)$ computed for test image one.}
        \label{f:hypo0result}
  \end{subfigure}
	
  \begin{subfigure}[t]{.5\textwidth}
        \centering
        \includegraphics[width=.92\linewidth]{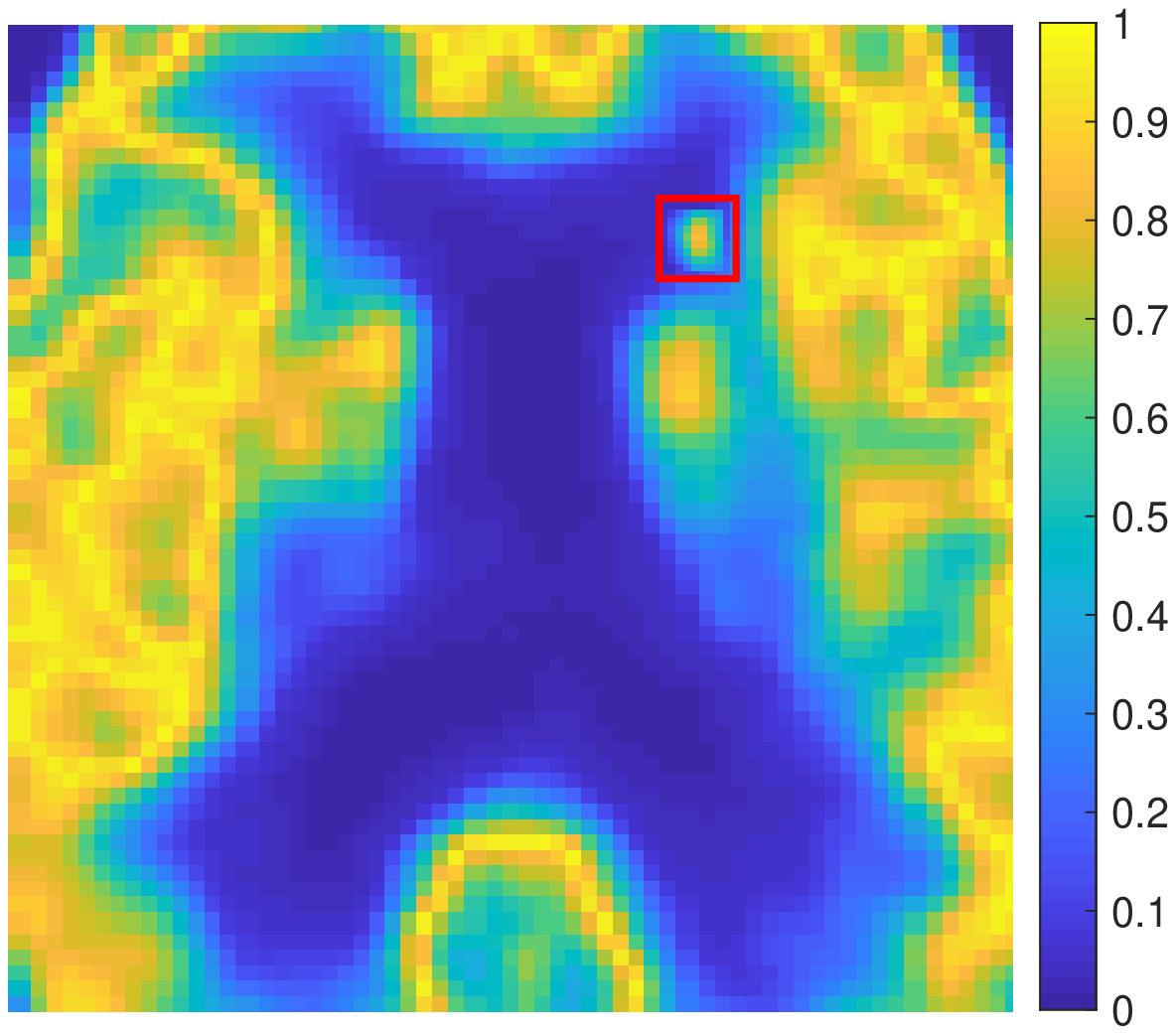}
\caption{Surrogate test image two.}
        \label{f:hypo1}
  \end{subfigure}
	~~~~
  \begin{subfigure}[t]{.5\textwidth}
        \centering
        \includegraphics[width=.92\linewidth]{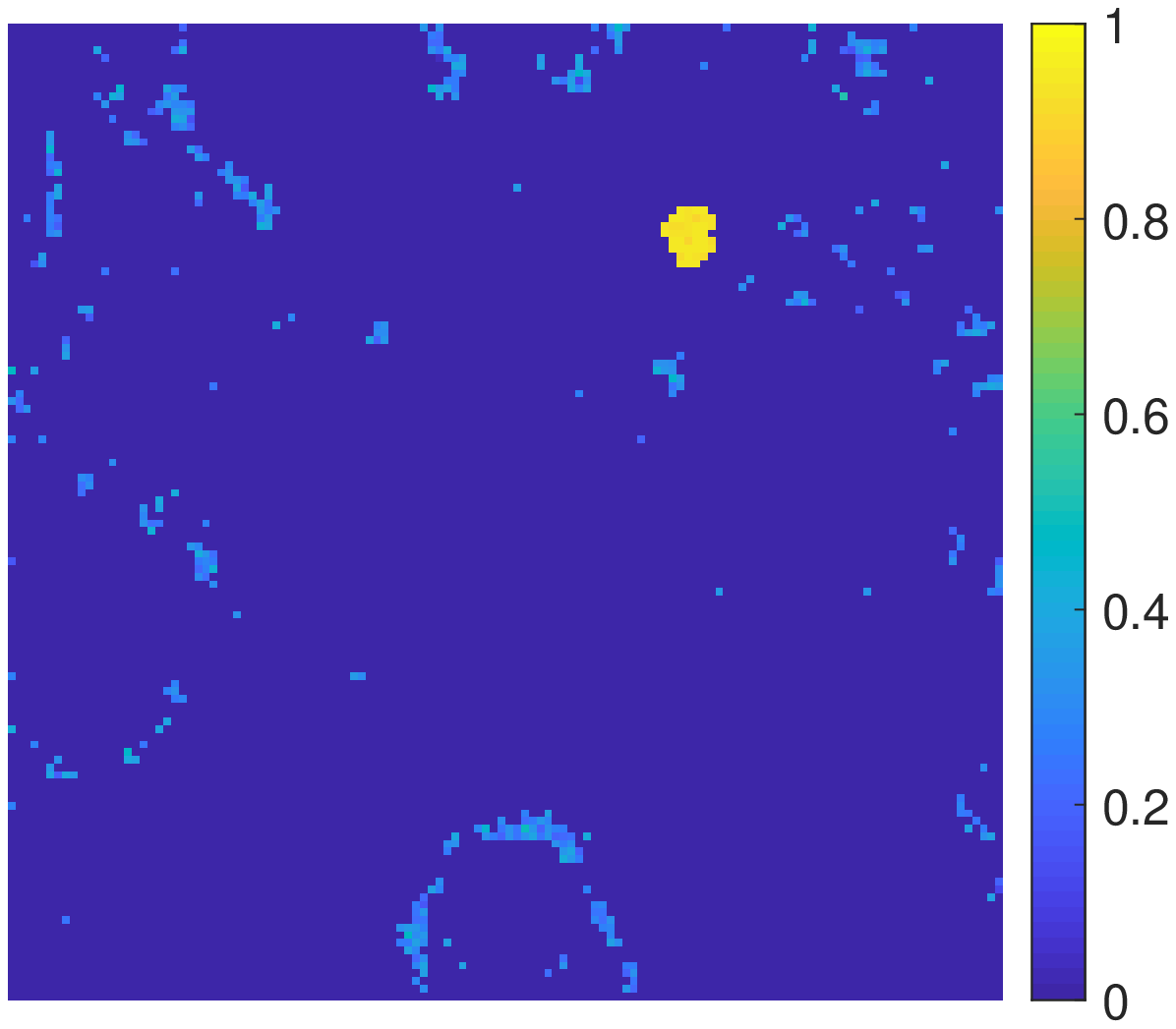}
\caption{The credible level $(1-\alpha)$ computed for test image two.}
          \label{f:hypo1result}
  \end{subfigure}

  \begin{subfigure}[t]{.5\textwidth}
        \centering
        \includegraphics[width=.92\linewidth]{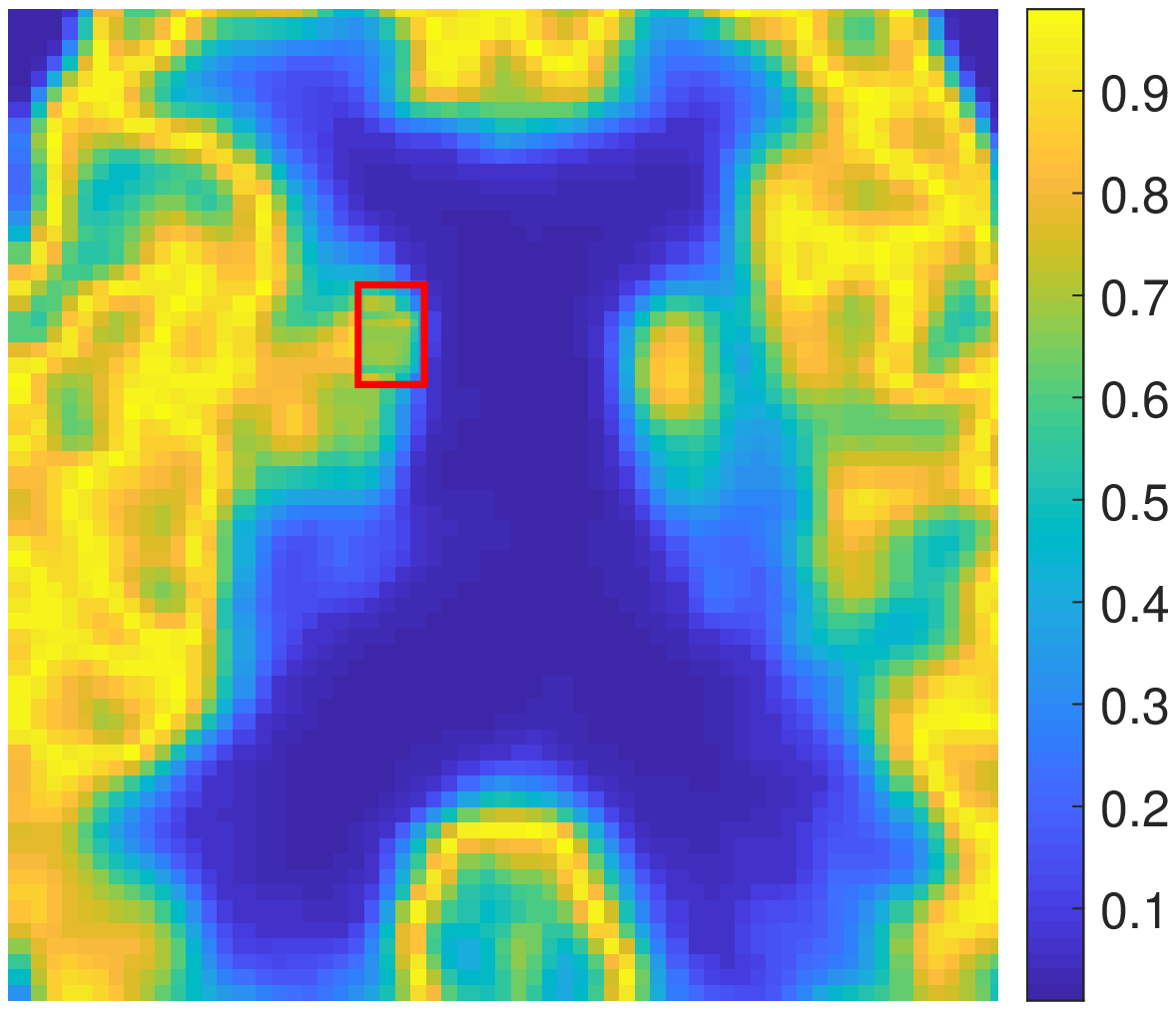}
\caption{Surrogate test image three.}
         \label{f:hypo2}
  \end{subfigure}
	~~~~
  \begin{subfigure}[t]{.5\textwidth}
    \centering
    \includegraphics[width=.92\linewidth]{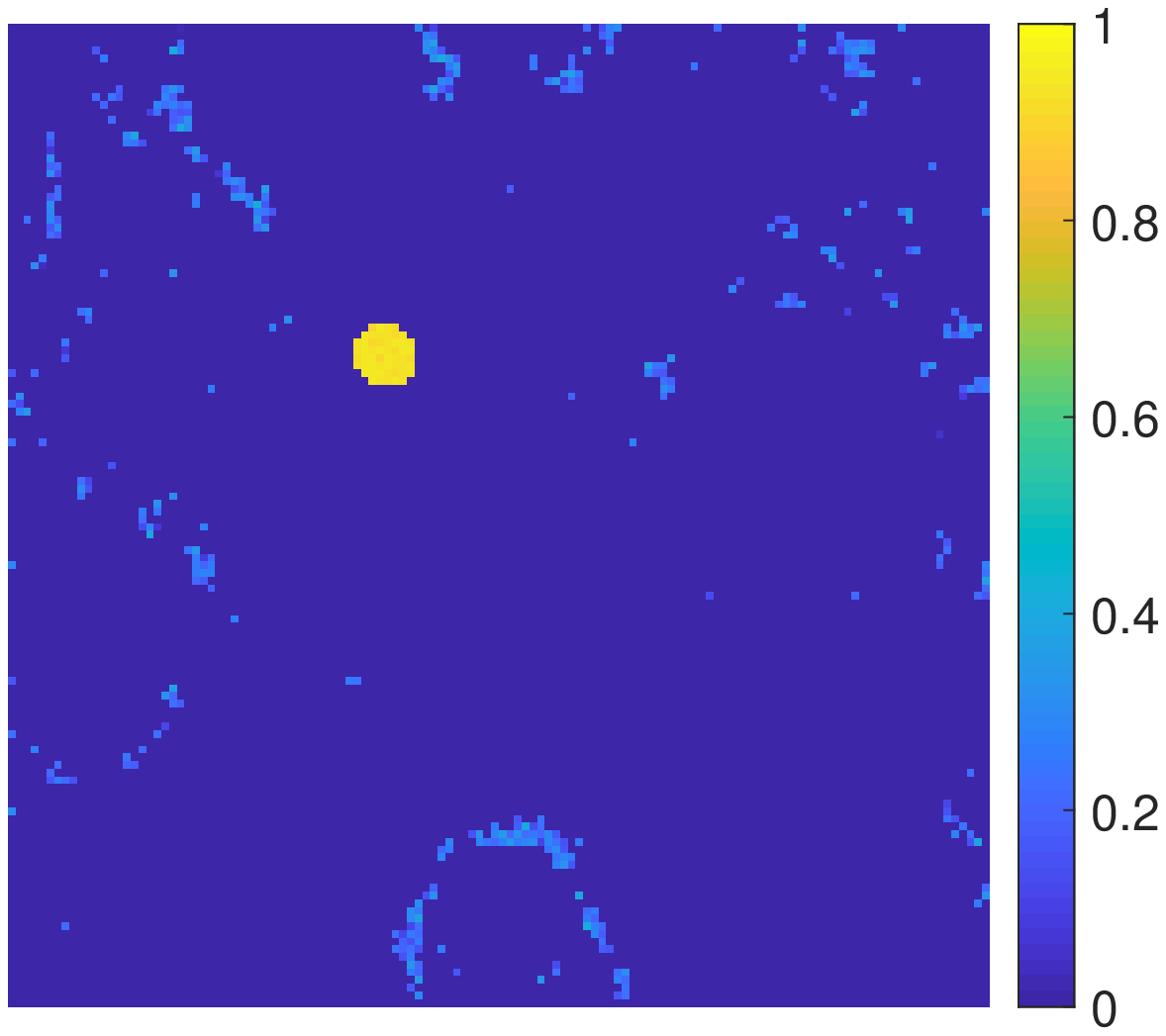}
\caption{The credible level $(1-\alpha)$ computed for test image three.}
     \label{f:hypo2result}
  \end{subfigure}
  \caption{($K=1$) The credible level for three surrogate test images.}
  \label{f:hysubadd}
\end{figure}

\begin{figure}[!htp]
  \begin{subfigure}[t]{.5\textwidth}
        \centering
        \includegraphics[width=.92\linewidth]{figs/hypo_noise_image}
\caption{Surrogate test image one.}
        \label{f:hypo0v2}
  \end{subfigure}
	~~~~
  \begin{subfigure}[t]{.5\textwidth}
        \centering
        \includegraphics[width=.92\linewidth]{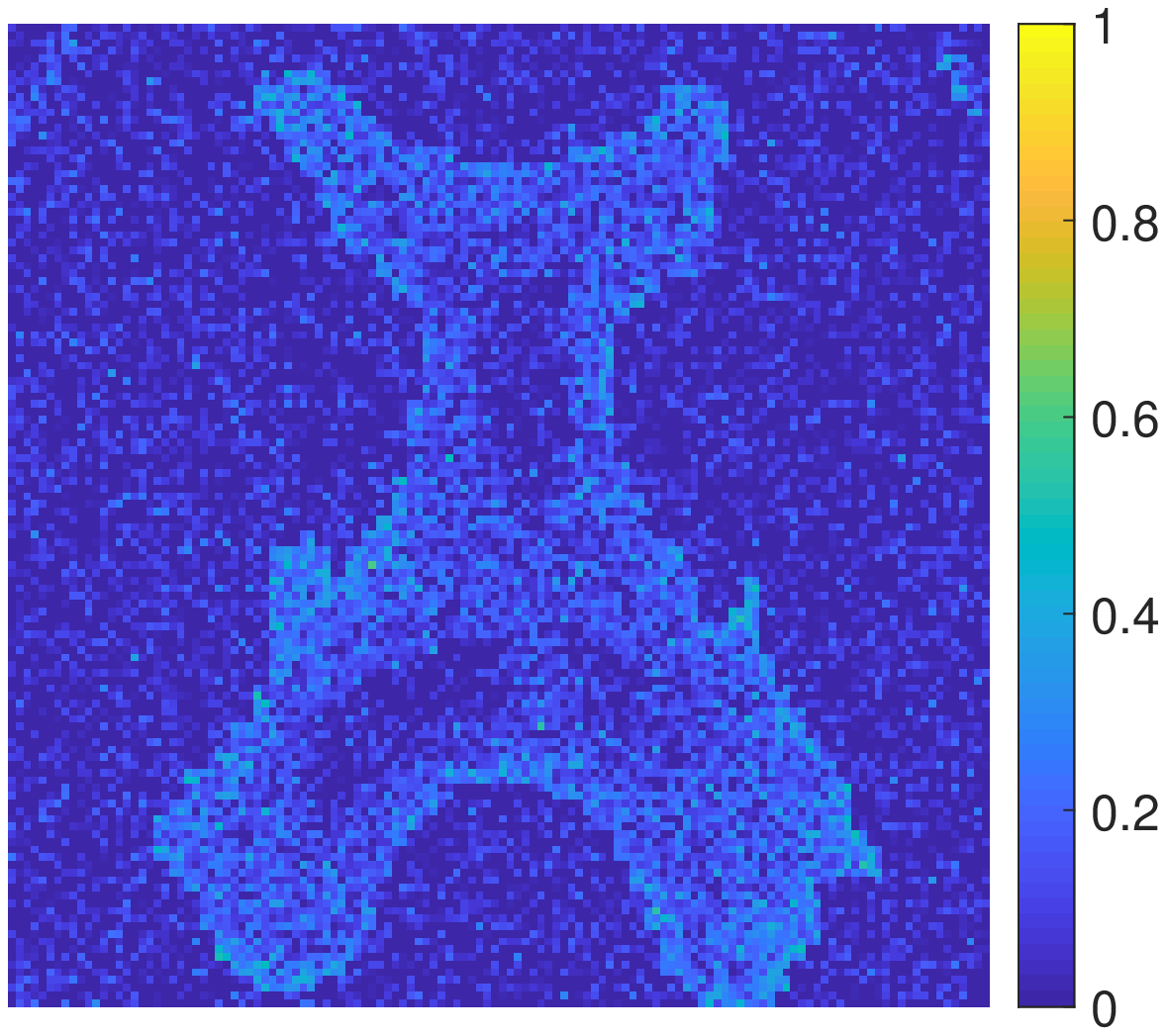}
\caption{The credible level $(1-\alpha)$ computed for test image one.}
        \label{f:hypo0resultv2}
  \end{subfigure}
	
  \begin{subfigure}[t]{.5\textwidth}
        \centering
        \includegraphics[width=.92\linewidth]{figs/hypo_fake_image_add}
\caption{Surrogate test image two.}
        \label{f:hypo1v2}
  \end{subfigure}
	~~~~
  \begin{subfigure}[t]{.5\textwidth}
        \centering
        \includegraphics[width=.92\linewidth]{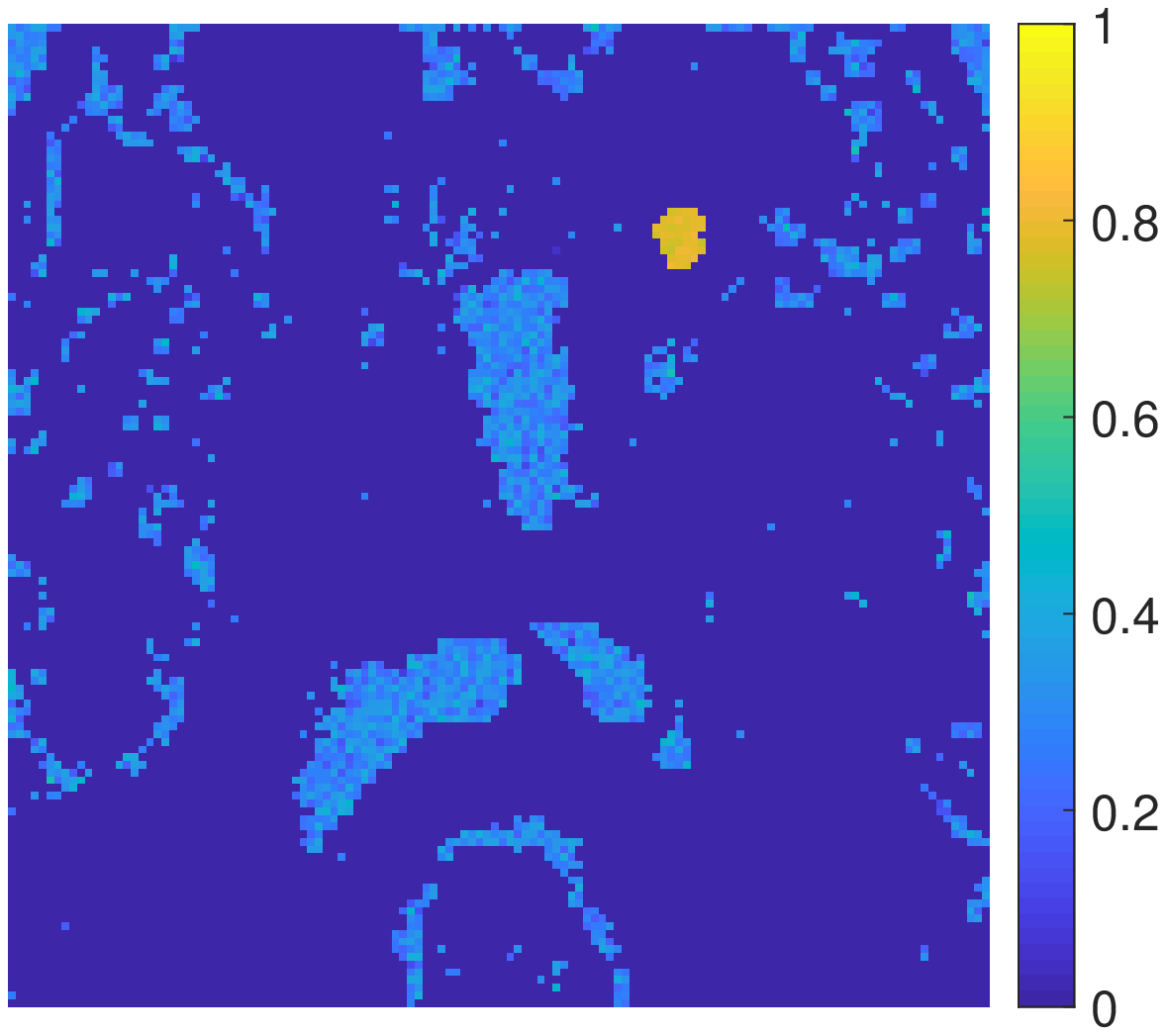}
\caption{The credible level $(1-\alpha)$ computed for test image two.}
          \label{f:hypo1resultv2}
  \end{subfigure}

  \begin{subfigure}[t]{.5\textwidth}
        \centering
        \includegraphics[width=.92\linewidth]{figs/hypo_fake_image_sub}
\caption{Surrogate test image three.}
         \label{f:hypo2v2}
  \end{subfigure}
	~~~~
  \begin{subfigure}[t]{.5\textwidth}
    \centering
    \includegraphics[width=.92\linewidth]{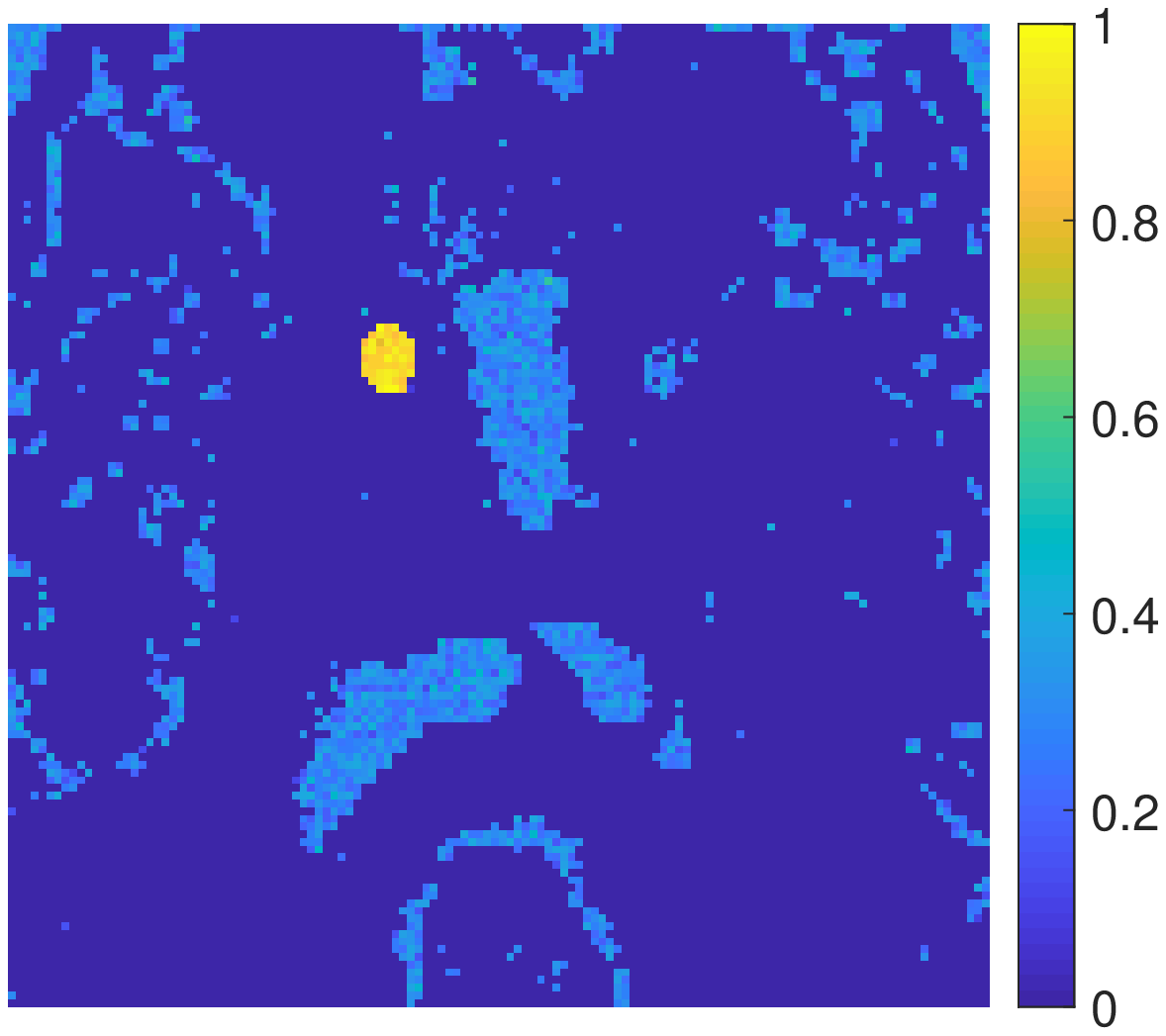}
\caption{The credible level $(1-\alpha)$ computed for test image three.}
     \label{f:hypo2resultv2}
  \end{subfigure}
  \caption{($K=0.5$) The credible level for three surrogate test images.}
  \label{f:hysubaddv2}
\end{figure}

\subsection{Identifying artifacts using HPDI}
As is discussed in Section~\ref{sec:artifact}, an important application of the proposed Bayesian framework is that the resulting posterior distribution 
can be used to detect artifacts in a reconstructed image. We now demonstrate this application with three surrogate test images
which are generated by making certain modification of the ground truth. 
Figs~\ref{f:hysubadd} summarize the results for $K=1$. 
Specifically the first image shown in Fig.~\ref{f:hypo0} is generated by adding some random noise to the ground truth without 
any structural changes, the second one shown in Fig.~\ref{f:hypo1} is generated by
 adding some artificial components to the ground truth, and the third one shown
in Fig.~\ref{f:hypo2}, is the result of removing some components from the ground truth.
Thus both the last two test images have structural changes from the ground truth, 
and in both figures, the regions in which components are altered from the ground truth are highlighted with red boxes. 
We compute the credible level $(1-\alpha)$ for all three images, 
and show the results in Figs.~\ref{f:hypo0result} (for test image 1), \ref{f:hypo1result} (for test image 2) and \ref{f:hypo2result} (for test image 3).
 It can be seen here that, though the first test image is visibly perturbed by random noise, it does not have structural difference from the 
ground truth and so credible level result in Fig.~\ref{f:hypo0result} does not suggest any region has high likelihood to contain artifacts. 
On the other hand, in the other two test images, 
 at the locations where the original image is altered (i.e., artifacts introduced), the resulting credible level $(1-\alpha)$ is significantly higher 
than other regions, suggesting that these locations may contain artifacts. 
Figs.~\ref{f:hysubaddv2} shows the same test results but for $K=0.5$, and one can see that the figures exhibits qualitative the same behaviors as 
those of $K=1$.  
The results demonstrate that the proposed method can rather effectively detect the artifacts in a test image. 

\section{Conclusions}
In this work, we have presented a complete treatment for performing
Bayesian inference and uncertainty quantification for medical image reconstruction problems with Poisson data.
In particular, we formulate the problem in an infinite dimensional setting and
we prove that the resulting posterior distribution is well-posed in this setting.
Second, to sample the unknown function/image, we provide a modified pCNL MCMC algorithm,
the efficiency of which is independent of discretization dimensionality.
Specifically the modified algorithm calculates
the offset direction in the original pCNL algorithm by using a primal-dual method, to avoid computing the  gradient of the TV term in our formulation.
Third, we also give a method to determine the TV regularization parameter $\lambda$ which is critical for the prior distribution.
The method is based on the realized discrepancy method for assessing model fitness.
Finally we provide an application of the uncertainty information obtained by the Bayesian framework, using
the posterior distribution to identify possible artifacts in an image reconstructed.
We believe the proposed Bayesian framework can be used to reconstruct images and evaluate the uncertainty associated to reconstruction the results in many practical medical imaging problems with Poisson data.

There are several problems related to this work that we plan to investigate in the future.
In the future, we plan to apply the methods developed in this work to those real-world problems, especially the PET image reconstruction.

\appendix
\section{Proof of Proposition~\ref{prop:phi}}
We provide a proof of Proposition~\ref{prop:phi} here.
\begin{proof}
(1) From Eq.~\eqref{e:radon}, and Eq.~\eqref{e:ubound}, we obtain directly that
\[
0<\bm{\underline{\theta}}
\leq \bm{\theta}
\leq \bm{\overline{\theta}},
\]
 for two constant vectors $\underline{\bm{\theta}}$ and $\bm{\overline{\theta}}$.
 It follows directly that  
$ \|\ln \bm{\theta}\|_2\leq l_{\max}$ for a positive constant $l_{\max}$.

For every $r>0$,  we have  $\|y\|_2<r$.  By the Cauchy-Schwarz inequality, we obtain the lower bound,
\begin{align*}
   \Phi(z) &=\<\bm{\theta},\bm{1}\> - \< y,\ln\bm{\theta}\>  \notag \\
   &\geq  \<\bm{\theta},\bm{1}\>  -\| \ln\bm{\theta}\|_2 \|y\|_2 \notag \\
   &\geq \<\bm{\underline{\theta}},\bm{1}\> - l_{\max}\|y\|_2  \notag \\
   & \geq \<\bm{\underline{\theta}},\bm{1}\> - l_{\max}r.
\end{align*}
For upper bound, once again we apply the Cauchy-Schwarz inequality to the functional $\Phi$, obtaining
\begin{align*}
\Phi(z) &=\<\bm{\theta},\bm{1}\> - \< y,\ln\bm{\theta}\> \\
&\leq  \<\bm{\theta},\bm{1}\>  +\| \ln\bm{\theta}\|_2 \|y\|_2 \notag \\
&\leq \<\bm{\overline{\theta}},\bm{1}\> + l_{\max}\|y\|_2  \notag \\
   & \leq \<\bm{\overline{\theta}},\bm{1}\> + l_{\max}r.
\end{align*}

(2) In this proof, we use $M$ for positive constants. Let $z$ and $v$ be any two elements in $X$,
and we have,
\begin{align} \label{e:phizv2}
|\Phi(z)-\Phi(v)| &=  |\< Af(z)-Af(v), \bm{1}\>  -\< y, \ln(Af(z))-\ln(Af(v)) \> |\notag\\ 
&\leq | \< Af(z)-Af(v), \bm{1}\> | + |\< y, \ln(Af(z))-\ln(Af(v)) \> |\notag\\
&\leq \|\bm{1}\|_2 \|Af(z)-Af(v)\|_2 +\|y\|_2 \|\ln(Af(z))-\ln(Af(v))\|_2\notag\\
&\leq \|\sqrt{d} \|Af(z)-Af(v)\|_2 +\|y\|_2 M \|Af(z)-Af(v)\|_2\notag\\
&=  (\sqrt{d} +\|y\|_2 M) \|Af(z)-Af(v)\|_2.
\end{align}
Since the Radon transform $A$ is a bounded linear operator from the $L_2$ space to $R^d$~\cite{natterer2001mathematics}, we have,  
\begin{align*}\label{e:af2}
\|Af(z)-Af(v)\|_2
 & \leq \|A\| \|f(z)-f(v)\|_{L_2}\\
 & =\|A\| \|\int^z_v e^{-t^2}dt \|_{L_2}\leq \|A\|\| \int ^z_v dt\|_{L_2} = \|A\|\|z-v\|_{L_2},
\end{align*}
which completes the proof.

(3) For any $y,\,y'\in Y$, it is easy to show,
\[
|\Phi(z,y)-\Phi(z,y')|= | \< y-y', \ln \bm{\theta} \> | \leq \|y-y'\|_2\|\ln\bm{\theta}\|_2\leq l_{\max}\|y-y'\|_2.\]

\end{proof}

\section{Proof of Theorem \ref{thm:db}}
We define $\eta_0(z,v)$ to be the measure $\eta(z,v)$ on $X\times X$ with $\Psi\equiv 0$,
and it is obvious that the measure $\eta_0(z,v)$ is Gaussian.
Moreover we have, 
\begin{equation}
\eta(dz,dv)=q(z,dv)\mu(dz),\quad
\eta_0(dz,dv)=q(z,dv)\mu_0(dz),
\end{equation}
and that the measures $\mu$ and $\mu_0$ are equivalent. It follows that $\eta$ and $\eta_0$ are equivalent and 
\begin{equation}\label{eq:priu}
\frac{d\eta}{d\eta_0}(z,v)=\frac{d\mu}{d\mu_0}(z)=Z\exp(-\Psi(z)),\quad
\frac{d\eta}{d\eta_0}(v,z)=Z\exp(-\Psi(v)).
\end{equation}
Now we define
\[\eta_0^{\perp}(dz,dv)=q(v,dz)\mu_0(dv),\]
and by some elementary calculations we can derive, 
\begin{equation}\label{eq:like}
\begin{aligned}
\frac{d\eta_0^{\perp}}{d\eta_0}
(z,v)
=&\exp(-\frac{1}{2}\frac{||2\mathcal{C}_0^{-\frac{1}{2}}(z-v)+\delta\mathcal{C}_0^{-\frac{1}{2}}(v+z)+2\delta \mathcal{C}_0^{\frac{1}{2}}g(v)||^2}{8\delta}-\frac{1}{2}\frac{||v||^2}{\mathcal{C}}\\
&\ \ \ \ \ \ \ \ \ \ \ \ +\frac{1}{2}\frac{||2\mathcal{C}_0^{-\frac{1}{2}}(v-z)+\delta\mathcal{C}_0^{-\frac{1}{2}}(z+v)+2\delta \mathcal{C}_0^{\frac{1}{2}}g(z)||^2}{8\delta}+\frac{1}{2}\frac{||z||^2}{\mathcal{C}})\\
=&\exp(-\frac{1}{2}\langle z-v,g(v)\rangle-\frac{\delta}{4}\langle(z+v),g(v)\rangle-\frac{\delta}{4}\langle \mathcal{C}_0^{\frac{1}{2}}g(v),\mathcal{C}_0^{\frac{1}{2}}g(v)\rangle\\
&\ \ \ \ \ \ \ \ \ \ \ \ +\frac{1}{2}\langle v-z,g(z)\rangle+\frac{\delta}{4}\langle(v+z),g(z)\rangle+\frac{\delta}{4}\langle \mathcal{C}_0^{\frac{1}{2}}g(z),\mathcal{C}_0^{\frac{1}{2}}g(z)\rangle).
\end{aligned}
\end{equation}
As $z,v\in X$ and $g(z)$ is in the Cameron-Martin space of $\mu_0$, $\langle z-v,g(z)\rangle$, $\langle v+z,g(z)\rangle$ and $||\mathcal{C}_0^{\frac{1}{2}}g(z)||^2$ are finite,
and $\frac{d\eta_0^{\perp}}{d\eta_0}$ is well defined. Now recall that, 
\begin{equation}\label{eq:1}
\begin{aligned}
\frac{d\eta^{\perp}}{d\eta}(z,v)=\frac{d\eta}{d\eta_0}(v,z)\frac{d\eta^{\perp}_0}{d\eta_0}(z,v)\frac{d\eta_0}{d\eta}(z,v).
\end{aligned}
\end{equation}
Substituting Eqs.~\eqref{eq:priu} and \eqref{eq:like} into the Eq. \eqref{eq:1} yields, 
\begin{equation}
\frac{d\eta^{\perp}}{d\eta}(z,v)=\exp(\rho(z,v)-\rho(v,z)),
\end{equation}
where
\begin{equation}
\rho(z,v)=\Phi(z)+\frac{1}{2}\langle v-z,g(z)\rangle+\frac{\delta}{4}\langle z+v,g(z)\rangle+\frac{\delta}{4}||\mathcal{C}_0^{\frac{1}{2}}g(z)||^2.
\end{equation}
\bibliographystyle{plain}
\bibliography{pet}
\end{document}